\documentclass[notitlepage]{article}

 \usepackage[matrix,arrow,curve]{xy}
 \usepackage{graphicx}
 \usepackage[english]{babel}
 \usepackage{amsmath}
 \usepackage{amssymb, indentfirst}
 \usepackage{amsthm}
 \usepackage{enumitem}
\begin{document}

\newtheorem{cor}{Corrolary}
\newtheorem{lem}{Lemma}
\newtheorem{rem}{Remark}
\newtheorem{ex}{Example}
\newtheorem{opr}{Definition}
\newtheorem{thm}{Theorem}
\newtheorem{prop}{Proposition}

\title{On the derived Picard group of the Brauer star algebra}
\author{Alexandra Zvonareva \let\thefootnote\relax\footnote{This work was partially supported by RFFI 13-01-00902 and by the Chebyshev Laboratory  (Department of Mathematics and Mechanics, St. Petersburg State University) under RF Government grant 11.G34.31.0026.}}
\date{}
\maketitle

\begin{abstract}
In this paper we show that the derived Picard group ${\rm TrPic}(A)$
of the Brauer star algebra of type $(n,t)$ is generated by shift,
${\rm Pic}(A)$ and equivalences $\{H_i\}_{i=1}^n$ in the case $t>1$,
where $H_i$ were shown to satisfy the relations of the braid group
on the affine diagram $\widetilde{A}_{n-1}$ by Schaps and
Zakay-Illouz. In the multiplicity free case we show that ${\rm
TrPic}(A)$ is generated by a slightly bigger set.
\end{abstract}

\section{Introduction}

The derived Picard group ${\rm TrPic}(A)$ of an algebra $A$ is the
group of isomorphism classes of two-sided tilting complexes in
$D^b(A\otimes A^{op})$, with the product of the classes of $X$ and
$Y$ given by the class $X \otimes_{A} Y$. Equivalently ${\rm
TrPic}(A)$ is the group of the standard autoequivalences of $D^b(A)$
modulo natural isomorphisms. It is easy to see that the derived
Picard group is invariant under derived equivalence.

Rouquier and Zimmermann started the study of the derived Picard
group of Brauer tree algebras \cite{RZ}. In the case of multiplicity
one they constructed a morphism from Artin's braid group on $n+1$
strings (n is the number of simple $A$-modules) to ${\rm TrPic}(A)$
and showed it to be an isomorphism modulo some central subgroup when
$n=2$. In \cite{Zi} Zimmermann generalized these results to the case
of arbitrary multiplicity. Khovanov and Seidel defined an action of
the Artin's braid group on the bounded derived category of a certain
algebra similar to the Brauer tree algebra, according to their
results the action of the Artin's braid group on the bounded derived
category of a Brauer tree algebra with multiplicity one is faithful
\cite{KS}.

Schaps and Zakay-Illouz constructed an action of the braid group on
the affine diagram $\widetilde{A}_{n-1}$ on the bounded derived
category of a Brauer tree algebra with arbitrary multiplicity
\cite{SI}. Muchtadi-Alamsyah showed this action to be faithful in
the case of multiplicity one \cite{IM}.

Schaps and Zakay-Illouz also raised the question whether in the case
of multiplicity $\neq 1$ the braid generators, together with the
shift and ${\rm Pic}(A)$, generate the entire derived Picard group
and whether the homomorphism from the braid group is one-to-one.
Using the technique of tilting mutations developed by Aihara and
Iyama \cite{AI}, \cite{Ai2}, we answer the first question
positively.

Let $A$ be the Brauer star algebra of
type (n,t), i.e. $A$ is the Brauer tree algebra corresponding to a star with $n$ edges and the exceptional vertex in the middle, the multiplicity of the exceptional vertex is $t$ (for more details see section \ref{dereq}). Let $A= \oplus_{i \in \mathbb{Z}/n\mathbb{Z}} P_i$ be the decomposition of $A$ into indecomposable mutually non-isomorphic
projective modules. Let $\{H_i\}_{i=1}^n$ be standard autoequivalences of of $D^b(A)$ such that
$$H_i(P_j)=\left\{
            \begin{array}{ll}
              0\mbox{ } \rightarrow \mbox{ } 0\mbox{ } \mbox{ }\rightarrow \mbox{ } P_j, & j \neq i, i-1  \\
              0\mbox{ }\rightarrow \mbox{ } 0\mbox{ } \mbox{ }\rightarrow \mbox{ } P_{i-1}, & j = i \\
              P_i \xrightarrow{\beta} P_{i-1} \xrightarrow{soc} P_{i-1}, & j =
              i-1,
            \end{array}
          \right.
$$
where the right most non-zero terms of the complexes above are concentrated in degree $0$.

\textbf{Theorem 3.} \emph{ Let $A$ be the Brauer star algebra of
type (n,t), $t>1$.  Then ${\rm TrPic}(A)$ is generated by shift,
${\rm Pic}(A)$ and equivalences $H_i$.}

In the multiplicity free case we show that ${\rm TrPic}(A)$ is
generated by a slightly bigger set.

\textbf{Acknowledgement:} I would like to thank Rapha\"el Rouquier
for bringing my attention to tilting mutations and Mikhail Antipov,
Yury Volkov and Alexander Generalov for numerous discussions and for
valuable suggestions on the exposition.

\section{Preliminaries}
\subsection{Derived equivalences}\label{dereq}

Let $A$ and $B$ be algebras over a commutative ring $R$ and let $A$
and $B$ be projective as modules over $R$. By $A$-modules we mean
left $A$-modules. In the path algebra of a quiver the product
of arrows $\overset{a}{\rightarrow} \overset{b}{\rightarrow}$ will
be denoted by $ab.$ Denote by $C(A)$ the category of
complexes of $A$-modules, by $D^b(A)$ the bounded derived category of
$A$, by $K^b({\rm proj}-A)$ the homotopy category of bounded complexes of
finitely generated projective modules. For an object $T$ of some additive category $\mathcal{T}$ denote by ${\rm add}(T)$ the smallest full subcategory of $\mathcal{T}$ which is closed under finite direct sums, summands and isomorphisms and contains $T$. The following theorem gives a necessary and sufficient condition for $A$ and $B$ to be derived equivalent.

\begin{thm} \emph{(Rickard, Keller, \cite{Ri1}, \cite{Ri3},
\cite{Ke})} The following are equivalent: \begin{enumerate}
                  \item The categories $D^b(A)$ and $D^b(B)$ are
                  equivalent as triangulated categories.
                  \item The categories $K^b({\rm proj}-A)$ and $K^b({\rm proj}-B)$ are
                  equivalent as triangulated categories.
                  \item There is a complex $T \in K^b({\rm proj}-A)$ such that \begin{itemize}
                                                                     \item
                                                                     ${\rm Hom}_{D^b(A)}(T,
                                                                     T[i])=0$
                                                                     for
                                                                     $i \neq
                                                                     0$,
                                                                     \item
                                                                     $K^b({\rm proj}-A)$
                                                                     is
                                                                     generated
                                                                     by
                                                                     ${\rm add}(T)$ as a triangulated category,
                                                                     \item
                                                                     ${\rm End}_{D^b(A)}(T)\simeq
                                                                     B^{op}$.
                                                                   \end{itemize}
                  \item There is a bounded complex $X$ of $(A\otimes
                  B^{op})$-modules whose restrictions to $A$ and to $B^{op}$ are perfect and a bounded complex $Y$ of $(B\otimes
                  A^{op})$-modules whose restrictions to $B$ and to $A^{op}$ are
                  perfect such that $X\otimes_{B}Y \simeq A$ in
                  $D^b(A\otimes A^{op})$ and $Y\otimes_{A}X \simeq B$ in
                  $D^b(B\otimes B^{op})$.
                \end{enumerate}\end{thm}
The complex $T$ from (3) is called a tilting complex, $X$ and $Y$
from (4) are called two-sided tilting complexes inverse to each
other. The complex $X$ viewed as a complex of $A$-modules or as a
complex of $B^{op}$-modules is a tilting complex. The inverse
equivalences between $D^b(A)$ and $D^b(B)$ are given by
$X\otimes_{B}-$ and $Y\otimes_{A}-$. Such equivalences are called
standard.

Let $T$ be a tilting complex such that ${\rm End}_{D^b(A)}(T)\simeq B^{op}$. There
exists a two-sided tilting complex $X$ of $(A\otimes
B^{op})$-modules whose restriction to $A$ is isomorphic to $T$ in
$D^b(A)$. If $X'$ is another two-sided tilting complex of $(A\otimes
B^{op})$-modules whose restriction to $A$ is isomorphic to $T$, then
there exists $\sigma \in {\rm Aut}(B)$ such that $X'= X \otimes_{B}
B_{\sigma} $, where $B_{\sigma}$ is a $(B\otimes B^{op})$-module
isomorphic to $B$ but with the right action twisted by $\sigma$
\cite{RZ}.

\begin{opr}

Let $\Gamma$ be a tree with $n$ edges and a distinguished vertex,
which has an assigned multiplicity $t \in \mathbb{N}$ (this vertex
is called exceptional, $t$ is called the multiplicity of the
exceptional vertex). Let us fix a cyclic ordering of the edges
adjacent to each vertex in $\Gamma$ (if $\Gamma$ is embedded into
the plane, we will assume that the cyclic ordering is clockwise). In
this case $\Gamma$ is called a Brauer tree of type $(n,t)$.

\end{opr}

The case when the tree is a star and the exceptional vertex is in
the middle is called the Brauer star.

To a Brauer tree of type $(n,t)$ one can associate an algebra
$A(n,t)$. The algebra $A(n,t)$ is a path algebra of a quiver with
relations. Let us construct a Brauer quiver $Q_{\Gamma}$ using the
Brauer tree $\Gamma$. The vertices of $Q_{\Gamma}$ are the edges of
$\Gamma$. Let $i$ and $j$ be two edges incident to the same vertex
$x$ in $\Gamma$, denote by $\geq$ the cyclic order of the edges
incident to $x$. If for any $h$, incident to $x$, such that $j \geq
h \geq i$ either $h=j$ or $h=i$, then there is an arrow from the
vertex $i$ to the vertex $j$ in $Q_{\Gamma}$. $Q_{\Gamma}$ has the
following properties: $Q_{\Gamma}$ is the union of oriented cycles
corresponding to the vertices of $\Gamma$, each vertex of
$Q_{\Gamma}$ belongs to exactly two cycles. The cycle corresponding
to the exceptional vertex is called exceptional. The arrows of
$Q_{\Gamma}$ can be divided into two families $\alpha$ and $\beta$
in such a manner that the arrows belonging to intersecting cycles
are in different families. By a slight abuse of notation arrows
belonging to the family $\alpha$ will be denoted by $\alpha$ and
arrows belonging to the family $\beta$ will be denoted by $\beta$
respectively.

\begin{opr}
Let $k$ be an algebraically closed field. The basic Brauer tree
algebra $A(n,t)$, corresponding to a tree $\Gamma$ of type $(n,t)$
is isomorphic to $kQ_{\Gamma}/I$, where the ideal $I$ is generated
by the relations:
\begin{enumerate}
  \item $\alpha\beta=0=\beta\alpha$ for any arrows belonging to families $\alpha$ and $\beta$ respectively;
  \item for any vertex
  $i$, not belonging to the exceptional cycle,
  $e_i\alpha^{i_{\alpha}}e_i=e_i\beta^{i_{\beta}}e_i$, where $i_{\alpha},$ resp.
  $i_{\beta}$ is the length of the $\alpha$, resp. $\beta$-cycle,
  containing
  $i$ and $e_i$ is the idempotent corresponding to $i$;
  \item for any vertex
  $i$, belonging to the exceptional $\alpha$-cycle (resp. $\beta$-cycle),
  $e_i(\alpha^{i_{\alpha}})^te_i=e_i\beta^{i_{\beta}}e_i$ (resp.
  $e_i\alpha^{i_{\alpha}}e_i=e_i(\beta^{i_{\beta}})^te_i$).
\end{enumerate}

An algebra is called a Brauer tree algebra of type $(n,t)$ if it is
Morita equivalent to the algebra $A(n,t)$ for some Brauer tree
$\Gamma$ of type $(n,t)$.
\end{opr}

Note that the ideal $I$ is not admissible. Let $x$ be a leaf in $\Gamma$ which is not an exceptional vertex, let $i$ be the edge incident to $x$. The loop corresponding to $x$ in $Q_{\Gamma}$ is equal to the other cycle passing through the vertex corresponding to $i$, we will not draw such superfluous loops. From now on for
convenience all algebras are supposed to be basic. We will assume
that $n > 1$, since if $n=1$ the corresponding Brauer tree algebra
is local and by results of Rouquier and Zimmermann \cite{RZ} the
derived Picard group is generated by shift and ${\rm Pic}(A)$, this
restriction has no effect on the statement of the results. Note that
the Brauer tree algebras are symmetric.

Rickard showed that two Brauer tree algebras corresponding to the
trees $\Gamma$ and $\Gamma'$ are derived equivalent if and only if
their types $(n,t)$ and $(n',t')$ coincide \cite{Ri2} and it follows
from the results of Gabriel and Riedtmann that this class is closed
under derived equivalence \cite{GR}.

Let $A$ be a Brauer star algebra with $n$ edges and multiplicity
$t$. The quiver of $A$ is of the form:

$$
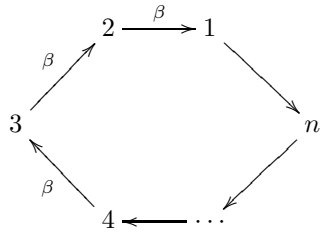
\xymatrix {
& 2 \ar[r]^-{\beta} & 1 \ar[dr] &  \\
    3 \ar[ur]^-{\beta} &  &  & n \ar[dl] \\
     & 4 \ar[ul]^-{\beta} & \cdots \ar[l] & \\
}
$$

Let us denote by $\{H_i\}_{i=1}^n$ some standard
autoequivalences of $D^b(A)$ which
act on the projective modules as follows:
\begin{equation}\label{**}
    H_i(P_j)=\left\{
            \begin{array}{ll}
              0\mbox{ } \rightarrow \mbox{ } 0\mbox{ } \mbox{ }\rightarrow \mbox{ } P_j, & j \neq i, i-1  \\
              0\mbox{ }\rightarrow \mbox{ } 0\mbox{ } \mbox{ }\rightarrow \mbox{ } P_{i-1}, & j = i \\
              P_i \xrightarrow{\beta} P_{i-1} \xrightarrow{soc} P_{i-1}, & j =
              i-1,
            \end{array}
          \right.
\end{equation}
where the right most non-zero terms of the complexes above are concentrated in degree $0$ and $soc$ is the morphism whose image is isomorphic to the socle
of $P_{i-1}$.

Schaps and Zakay-Illouz studied the
subgroup of the derived Picard group generated by $\{H_i\}_{i=1}^n$ \cite{SI2} and showed that $\{H_i\}_{i=1}^n$
satisfy the relations of the braid group on the affine diagram
$\widetilde{A}_{n-1}$.

\subsection{Mutations}\label{mut}

Let $k$ be an algebraically closed field. Let $\mathcal{T}$ be a
Krull-Schmidt, $k$-linear, Hom-finite triangulated category. A
morphism $X \overset{f}{\rightarrow} M' \in \mathcal{T} $ is called
left minimal if any morphism $g: M' \rightarrow M'$ satisfying $gf =
f$ is an isomorphism. Let $\mathcal{M}$ be a subcategory of
$\mathcal{T}$, $X$ an object of $\mathcal{T}$, $M'$ an object of
$\mathcal{M}$, a morphism $X \overset{f}{\rightarrow} M'$ is called
a left approximation of $X$ with respect to $\mathcal{M}$ if ${\rm
Hom}_{\mathcal{T}}(M',M)\overset{f^*}{\rightarrow}{\rm
Hom}_{\mathcal{T}}(X,M)$ is surjective for any $M \in \mathcal{M}$.
Right minimal morphisms and right $\mathcal{M}$-approximations are
defined dually.

$T \in \mathcal{T}$ is called silting if ${\rm
Hom}_{\mathcal{T}}(T,T[i])=0$ for any $i>0$ and $\mathcal{T}$ is
generated by ${\rm add}(T)$ as a triangulated category. We say that
a silting object $T$ is basic if $T$ is isomorphic to a direct sum of
indecomposable objects which are mutually non-isomorphic.

Let $T$ be a basic silting object in $\mathcal{T}$, $T=M\oplus X,$
$\mathcal{M}={\rm add}(M).$ Consider a triangle
\begin{equation}\label{*}
    X \overset{f}{\longrightarrow} M'
\longrightarrow Y \longrightarrow,
\end{equation}
where $f$ is a minimal left approximation of $X$ with respect to
$\mathcal{M}$. The morphism $f$ is unique up to isomorphism, since
if there is some other morphism $f'$ which is a minimal left
approximation of $X$ with respect to $\mathcal{M}$, then $\exists h$
such that $hf=f'$ and $\exists h'$ such that $h'f'=f$, hence
$f'=hh'f'$ and $f=h'hf$, but then by minimality $hh'$ and $h'h$ are
isomorphisms. The object $\mu^+_X(T):=M \oplus Y$ is called a left
mutation of $T$ with respect to $X$. Right mutations are defined
dually and are denoted by $\mu^-_X(T)$. By results of Aihara and
Iyama \cite{AI} $\mu^+_X(T)$ is again a basic silting object and
$\mu^-_Y(\mu^+_X(T))=T$ in the notation of the triangle (2). If $X$
is indecomposable the mutation is called irreducible.

Let $T,U$ be basic silting objects in $\mathcal{T}$. Set $T\geq U$
if ${\rm Hom}_{\mathcal{T}}(T,U[i])=0$ for any $i>0$. The relation
$\geq$ gives a partial order on the set of the isomorphism classes
of basic silting object of $\mathcal{T}$ \cite{AI}. We say that $U$
is connected (left-connected) to $T$ if $U$ can be obtained from $T$
by iterated irreducible (left) mutation. A triangulated category
$\mathcal{T}$ is called silting-connected if all basic silting
objects in $\mathcal{T}$ are connected to each other. $\mathcal{T}$
is strongly silting-connected if for any silting objects $T,U$ such
that $T\geq U$ the object $U$ is left-connected to $T$. Since in the
case of a symmetric algebra ${\rm Hom}_{D^b(A)}(T,T[i])\simeq D{\rm
Hom}_{D^b(A)}(T[i],T)$, where $D$ is the duality with respect to $k$
\cite{Ha}, then any silting object in $K^b({\rm proj}-A)$ is a
tilting complex, in this case instead of the term silting-connected
we will use the term tilting-connected.

\begin{thm} \emph{(Aihara, \cite{Ai2})} $K^b({\rm proj}-A)$ is
tilting-connected if $A$ is a representation-finite symmetric
algebra.
\end{thm}

Note that it also follows from \cite{Ai2} (Theorem 5.6 and Corollary
3.9) that in the case of representation-finite symmetric algebra the category
$K^b({\rm proj}-A)$ is strongly tilting connected. So any tilting
complex concentrated in non-positive degrees can be obtained from $A$
by iterated irreducible left mutations.

A well known example of tilting mutations is the mutations of Brauer
graphs or equivalently symmetric special biserial algebras. We will define it in the
context of Brauer trees, i.e. we will not consider loops.

Consider a Brauer tree algebra $A$ as a tilting complex over itself.
$A=(\bigoplus_{i=1, i \neq j}^{n}P_i) \oplus P_j $. First, assume
that $j$ is incident to a leaf. Consider a left mutation
$$\mu^+_{P_j}(A)=(\bigoplus_{i=1, i \neq j}^{n}P_i) \oplus (P_j
\overset{\beta}{\longrightarrow} P_l) $$
of $A$ with respect
to $P_j$ (the corresponding approximation is taken with respect to
${\rm add}( \bigoplus_{i=1, i \neq j}^{n}P_i)$), where
$P_l$ is the projective modules corresponding to the edge in the
Brauer tree which follows $j$ in the cyclic ordering of the edges
incident to the same vertex and $\beta$ corresponds to the arrow from $j$ to
$l$. The Brauer tree of $A$ is on
the left-hand side and the Brauer tree of
${\rm End}_{D^b(A)}(\mu^+_{P_j}(A))^{op}$ is on the right-hand side, the edge corresponding
to $P_j \overset{\beta}{\longrightarrow} P_l$ will be
also denoted by $j$.

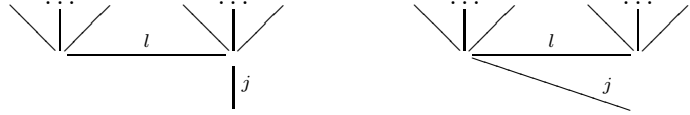
$$
\xymatrix @! @=0.35pc {
\ar@{-}[rd]&\dots \ar@{-}[d]&\ar@{-}[dl]&\ar@{-}[rd]&\dots \ar@{-}[d]&\ar@{-}[dl] && \ar@{-}[rd]&\dots \ar@{-}[d]&\ar@{-}[dl]&\ar@{-}[rd]&\dots \ar@{-}[d]&\ar@{-}[dl] \\
&\ar@{-}[rrr]^{l}&&&\ar@{-}[d]^{j}& && &\ar@{-}[rrr]^{l}\ar@{-}[drrr]^(.8){j}&&&&\\
&&&&& && &&&&& \\
}
$$

In the case where $j$ is not incident to a leaf in the Brauer tree of
$A$
$$\mu^+_{P_j}(A)=(\bigoplus_{i=1, i \neq j}^{n}P_i) \oplus (P_j
\overset{f}{\longrightarrow} P_m \oplus P_l), $$
where $P_m$ and
$P_l$ are the projective modules corresponding to the edges in the
Brauer tree which follow $j$ in the cyclic ordering of the edges
incident to the same vertices and $f={ \left(
            \begin {smallmatrix}
              \alpha \\
              \beta \\
            \end{smallmatrix}
          \right)
}$, where $\alpha$ and $\beta$ correspond to the arrows from $j$ to
$m$ and from $j$ to $l$ respectively. And on the level of Brauer trees we have:

$$
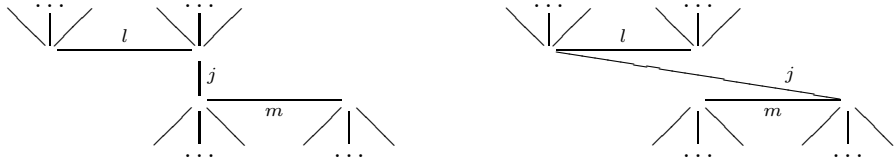
\xymatrix @! @=0.1pc {
\ar@{-}[rd]&\dots \ar@{-}[d]&\ar@{-}[dl]&\ar@{-}[rd]&\dots \ar@{-}[d]&\ar@{-}[dl]&&& && \ar@{-}[rd]&\dots \ar@{-}[d]&\ar@{-}[dl]&\ar@{-}[rd]&\dots \ar@{-}[d]&\ar@{-}[dl]&&& \\
&\ar@{-}[rrr]^{l}&&&\ar@{-}[d]^{j}&&&& && &\ar@{-}[rrr]^{l}\ar@{-}[rrrrrrd]^(.8){j}&&&&&&&\\
&&&&\ar@{-}[rd]\ar@{-}[d]\ar@{-}[ld]\ar@{-}[rrr]_{m}&&&\ar@{-}[rd]\ar@{-}[d]\ar@{-}[ld]& && &&&&\ar@{-}[rd]\ar@{-}[d]\ar@{-}[ld]\ar@{-}[rrr]_{m}&&&\ar@{-}[rd]\ar@{-}[d]\ar@{-}[ld]& \\
&&&&\dots&&&\dots& && &&&&\dots&&&\dots& \\
}
$$

As far as we know these moves were first introduced in \cite{KZ} and
\cite{Ka} but were also studied in \cite{Ant}, \cite{Ai1},
\cite{Du}, \cite{MS}. The mutations $\mu^-_{P_j}(A)$ are defined
dually and corresponds to the move in the opposite direction, i.e.
from the right-hand side to the left-hand side. The tilting
complexes of the form $\mu^{\pm}_{P_j}(A)$, where $A$ is an
arbitrary Brauer tree algebra, will sometimes be called elementary.
Mutations of an edge incident to a leaf will be called the mutations
of type \MakeUppercase{\romannumeral1} and mutations of an edge not
incident to a leaf will be called the mutations of type
\MakeUppercase{\romannumeral2}. Note also that these mutations
involve only edges: the exceptional vertex stays unchanged.

\section{Mutations and the derived Picard group}

Denote by ${\rm Pic}(A)$ the Picard group of an algebra $A$, i.e.
the group of isomorphism classes of invertible $A\otimes
A^{op}$-modules or equivalently the group of Morita autoequivalences
of $A$ modulo natural isomorphisms. The group ${\rm Out}(A)$ of outer autoequivalences of $A$
coincides with ${\rm Pic}(A)$ \cite{Bolla} and
is clearly a subgroup of ${\rm TrPic}(A)$. If $\sigma \in {\rm
Aut}(A)$, then there is an invertible $A\otimes A^{op}$-module
$A_{\sigma}$, where $A_{\sigma}$ is an $(A\otimes A^{op})$-module
isomorphic to $A$ as a left module but with the right action twisted
by $\sigma$. A bimodule $A_{\sigma}$ is isomorphic to $A_{\sigma'}$
if and only if $\sigma$ coincides with $\sigma'$ modulo the subgroup
of inner automorphisms. Consider equivalences $F, F': D^b(B)
\rightarrow D^{b}(A)$. Assume that they are given by two-sided
tilting complexes $X, X'$ whose restriction to $A$ is a tilting
complex $T$, then $X'= X \otimes_{B}
B_{\sigma}$ for some $\sigma \in {\rm Aut}(B)$ \cite{RZ}.

The proof of the following statement is trivial.

\begin{lem} Let $\mathcal{A}$, $\mathcal{B}$ be two triangulated
categories, $F: \mathcal{A} \rightarrow \mathcal{B}$ a triangular
equivalence, $\mathcal{M}$ a subcategory of $\mathcal{A}$, $X \in
\mathcal{A}$ and let $X \overset{f}{\rightarrow} M'$ ($M'
\overset{f}{\rightarrow} X$) be a minimal left (resp. right)
approximation of $X$ with respect to $\mathcal{M}$, then $F(X)
\overset{F(f)}{\longrightarrow} F(M')$ ($F(M')
\overset{F(f)}{\longrightarrow} F(X)$) is a minimal left (resp. right)
approximation of $F(X)$ with respect to $F(\mathcal{M})$.
\end{lem}

Let $A$ be a symmetric algebra and let $T$ be a basic tilting
complex over $A$. Denote by $B^{op}={\rm End}_{D^b(A)}(T)$ and by
$F$ a standard equivalence such that $F(B)\simeq T,$ where $F:
D^b(B) \rightarrow D^b(A).$ Let $T=M\oplus X,$ where $X$ is
indecomposable, let $\mu^{\pm}_X(T)$ be a tilting complex obtained
from $T$ by right or left mutation. Denote by $C^{op}={\rm
End}_{D^b(A)}(\mu^{\pm}_X(T))$ and by $G: D^b(C) \rightarrow D^b(A)$
a standard equivalence such that $G(C)\simeq \mu^{\pm}_X(T)$.
Consider $B$ as a tilting complex over itself. Summands of $T$ and
of $B$ are in one to one correspondence under $F$, we will denote by
$P_X$ the indecomposable projective $B$-module corresponding to $X,$
a tilting complex $\mu^{\pm}_{P_X}(B)$ is obtained from $B$ by
mutation with respect to $P_X$.

\begin{lem}In the above notation ${\rm End}_{D^b(B)}(\mu^{\pm}_{P_X}(B))\simeq C^{op}$.

\end{lem}

\textbf{Proof.} Consider only the case of right mutations.
$C^{op}={\rm End}_{D^b(A)}(M \oplus Y)$, ${\rm
End}_{D^b(B)}(\mu^{-}_{P_X}(B))={\rm End}_{D^b(B)}(F^{-1}(M) \oplus
Y')$, where $$Y' \rightarrow N \rightarrow P_X \rightarrow$$ is a
triangle, and $N \rightarrow P_X$ is a minimal right approximation
of $P_X$ with respect to ${\rm add}(F^{-1}(M))$. $F^{-1}$ is a
triangular equivalence, and due to the previous lemma $F^{-1}(Y)
\simeq Y'$, so ${\rm End}_{D^b(B)}(\mu^{-}_{P_X}(B))\simeq {\rm
End}_{D^b(B)}(F^{-1}(M \oplus Y))$ and the assertion
follows.\hfill\(\Box\)

Let's denote by $H:D^b(C) \rightarrow D^b(B)$ a standard equivalence
such that $H(C) \simeq \mu^{\pm}_{P_X}(B).$

\begin{lem} There exists $\sigma \in {\rm Aut}(C)$ such that
$G \simeq  F\circ H \circ (C_{\sigma} \otimes_{C}-)$.

\end{lem}

\textbf{Proof.} By the previous lemma we have the following diagram:
$$\xymatrix{D^b(C) \ar[dr]^G \ar[rr]^H&&D^b(B) \ar[dl]^F \\
&D^b(A)&}$$ We need to check that the action of $G$ and $F\circ H$
on $C$ coincides. We will deal only with the case of right mutation
again. $G(C)=\mu^{-}_X(T)$. Whereas, $H(C)=\mu^{-}_{P_X}(B) \simeq
F^{-1}(M) \oplus Y'$. Hence $F(H(C)) \simeq M \oplus F(Y')$, but
$F(Y') \simeq Y$, so $F(H(C)) \simeq M \oplus Y
=\mu^{-}_X(T)$.\hfill\(\Box\)

We will need the following technical lemma, the proof is
straightforward.

\begin{lem} Let $V=\cdots\rightarrow X'' \xrightarrow{x} X \xrightarrow{g} Y \xrightarrow{h} Z \rightarrow 0 \rightarrow \cdots$, $W=\cdots\rightarrow0 \rightarrow X' \xrightarrow{g'}
Y \xrightarrow{h'} Z'\xrightarrow{z} Z'' \rightarrow \cdots$ be two
objects of $D^b(A)$. Let
$$\xymatrix{\cdots  \ar[r] &X'' \ar[d] \ar[r]^{x} &X \ar[d]^{w} \ar[r]^{g} &Y \ar[d]^{{\rm Id}} \ar[r]^{h} & Z \ar[d]^{v} \ar[r] \ar[d] &0 \ar[r] \ar[d] &\cdots \\
\cdots \ar[r] &0 \ar[r] &X' \ar[r]^{g'} & Y \ar[r]^{h'} & Z'
\ar[r]^{z} &Z''\ar[r]&\cdots}$$ be a morphism $f \in D^b(A)$, then
$Cone(f) \simeq \cdots\rightarrow X'' \xrightarrow{x} X
\xrightarrow{w} X' \xrightarrow{hg'} Z \xrightarrow{v}
Z'\xrightarrow{z} Z'' \rightarrow \cdots $, where $Z'$ is concentrated
in the same degree as in $W$.
\end{lem}

\subsection{The standard construction of a tree}

In this section we will fix a standard way to build a tree from a
star using mutations. The procedure will involve only mutation of
type \MakeUppercase{\romannumeral1} (see section \ref{mut}).

Let $\Gamma$ be a Brauer tree of type $(n,t)$. Let us assume that
the root of $\Gamma$ is chosen in the exceptional vertex, and that
$\Gamma$ is embedded into the plane in such a manner that all non-root
vertices are situated on the plane lower than the root according to
their level (the further from the root, the lower, all vertices of
the same level lie on a horizontal line). The edges around vertices
are ordered clockwise.

Let $A$ be a Brauer star algebra. If the corresponding tree is
embedded into the plane as described above, let its edges be
labelled from left ro right. $A= \bigoplus_{i=1}^{n}P_i$, where
$P_i$ are indecomposable projective modules. We are going to perform
a series of irreducible mutations of $A$, after the $r$-th mutation
we are going to obtain a tilting complex $T^{r} =
\bigoplus_{i=1}^{n}T_{i}^{r}$. The summands of $T^{r}$ are labelled
as follows: $T^{0}=A$, $T_{i}^{0}=P_i$. Each mutation changes only
one summand of the tilting complex $T^{r-1}$, say $T_{j}^{r-1}$.
Denote by $T_{j}^{r}$ the summand which was changed by the $r$-th
mutation, $T_{i}^{r}:=T_{i}^{r-1}$, $i \neq j$. All this allows us
to write a composition of mutations. Denote by $\mu^{\pm}_j$ the
mutation which changes $T_{j}^{r-1}$. The Brauer tree of ${\rm
End}_{D^b(A)}(T^{r})^{op}$ can be obtained by the mutation of the
Brauer tree of $T^{r-1}$, for example by Lemma 2.

Note also, that since ${\rm End}_{D^b(A)}(T^{r})^{op}$ is a Brauer
tree algebra it is easy to compute the minimal left approximation of
$T_{j}^{r}$ with respect to other summands of $T^r$. If the edge $j$
corresponding to $T_{j}^{r}$ is not incident to a leaf in the Brauer
tree of ${\rm End}_{D^b(A)}(T^r)^{op}$, let $T_{m}^{r}$ and
$T_{l}^{r}$ be the summands corresponding to $m$ and $l$, the edges
following $j$ in the cyclic ordering and let $f$ be a morphism
corresponding to two arrows in ${\rm End}_{D^b(A)}(T^r)^{op}$, then
by Lemma 1 the morphism $T_{j}^{r} \overset{f}{\longrightarrow}
T_{m}^{r} \oplus T_{l}^{r}$ is a minimal left approximation of
$T_{j}^{r}$ with respect to other summands of $T^r$. If $j$ is
incident to a leaf, then $T_{j}^{r} \overset{f}{\longrightarrow}
T_{l}^{r}$ is a minimal left approximation of $T_{j}^{r}$, where $l$
is the only edge following $j$ in the cyclic ordering and $f$ is the
corresponding arrow. We will say that we mutate $j$ along the edges
$m$ and $l$ or along the the edge $l$ respectively. In addition if a
component of $f$ does not belong to the exceptional cycle, then any
non-zero morphism can be set as this component.

Let $\Gamma$ be a Brauer tree of type $(n,t)$, assume $n>1$. Let us
number the edges of the tree $\Gamma$ as follows: put 1 on the
left-hand edge incident to the root. If the edge with label 1 is not
incident to a leaf, put 2 on the edge incident to its
non-exceptional end which is the previous edge coming before the
edge with label 1 in the cyclic ordering; if the edge with label 1
is incident to a leaf, put 2 on the only edge which is the previous
edge coming before the edge with label 1 in the cyclic ordering.
Assume that the label $i$ is assigned to some edge, if it is not
incident to a leaf, put the label $i+1$ on the edge of the lower
level which is the previous edge coming before the edge with label
$i$ in the cyclic ordering (e.g.: edges $2$ and $3$ in Example
\ref{ex}); if the edge with the label $i$ is incident to a leaf, put
the label $i+1$ on the only edge which is the previous edge coming
before the edge with label $i$ in the cyclic ordering (e.g.: edges
$3$ and $4$ in Example \ref{ex}); if this edge has a label already,
find the edge $l$ with the biggest label which has an unlabelled
edge incident to its upper end $x$, put the label $i+1$ on the
previous edge coming before $l$ in the cyclic ordering around $x$
(e.g.: edges $4$ and $5$ in Example \ref{ex}). Note that the same
labelling can be obtained using Green walk or depth-first search
algorithm.

\begin{ex}\label{ex} The following labelling of the Brauer tree
$\Gamma$ is standard:

$$
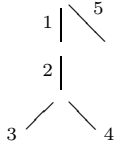
\xymatrix @! @=0.35pc {
      & \ar@{-}[dr]^{5} \ar@{-}[d]_{1} &   \\
      & \ar@{-}[d]_{2} &   \\
& \ar@{-}[dl] \ar@{-}[dr] &   \\
    ^3  &   & ^4\\
}$$

\end{ex}

Let $\phi_{\Gamma}: \{1,2,\dots,n\}\rightarrow \{0,1,\dots,n-1\}$ be a
function which assigns to a label $i$ the length of the shortest
path from this edge to the exceptional vertex, or equivalently the
level of the edge with label $i$ (assuming that the edges incident
to the exceptional vertex belong to the level 0).

\begin{lem}\label{l5} Let $\Gamma$ be a Brauer tree of type $(n,t)$ labelled as described above.
Let $T= (\mu^{+}_n)^{\phi_{\Gamma}(n)}
\circ(\mu^{+}_{n-1})^{\phi_{\Gamma}(n-1)}\circ\dots\circ
(\mu^{+}_1)^{\phi_{\Gamma}(1)}(A)$. Then the Brauer tree of ${\rm End}_{D^b(A)}(T)^{op}$ is $\Gamma$.

\end{lem}

\textbf{Proof.} Since $T= (\mu^{+}_n)^{\phi_{\Gamma}(n)}
\circ(\mu^{+}_{n-1})^{\phi_{\Gamma}(n-1)}\circ\dots\circ
(\mu^{+}_1)^{\phi_{\Gamma}(1)}(A)$, first we move edges with label $1$ as described in section \ref{mut}, then with label $2$ and so on. Thus we can build the tree of ${\rm End}_{D^b(A)}(T)^{op}$ by induction on the label. The desired description of the tree of ${\rm End}_{D^b(A)}(T)^{op}$ follows from the description of the moves and from the construction.\hfill\(\Box\)

\begin{ex} Let $T$ be the Brauer tree from Example \ref{ex}, $T= (\mu^{+}_4)^{2} \circ(\mu^{+}_{3})^{2}\circ \mu^{+}_2(A)$

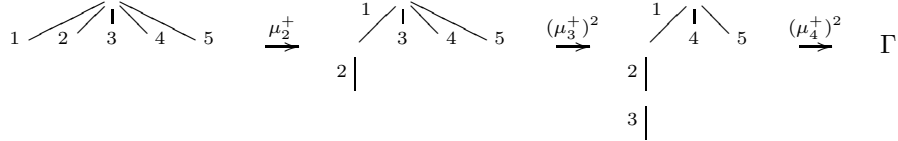
$$\xymatrix @! @=0.35pc{
      &   & \ar@{-}[dll] \ar@{-}[dl] \ar@{-}[d] \ar@{-}[dr] \ar@{-}[drr] & &  &   &   &   & \ar@{-}[dl]_1 \ar@{-}[d] \ar@{-}[dr] \ar@{-}[drr]  &   &   &   &   &   & \ar@{-}[dl]_1 \ar@{-}[d] \ar@{-}[dr] &   &   &      &    \\
    ^1 & ^2 & ^3 & ^4 & ^5 & \ar[r]^{\mu^{+}_2} & & \ar@{-}[d]_2  &  ^3 & ^4& ^5 & \ar[r]^{(\mu^{+}_3)^2} &  &  \ar@{-}[d]_2 & ^4  & ^5 &   \ar[r]^{(\mu^{+}_4)^2} &   &  \Gamma   \\
      &   &   &   &   &   &   &   &   &   &   & &  &  \ar@{-}[d]_3 &   &   &   &      &    \\
      &   &   &   &   &   &   &   &   &   &   &  & &   &   &   &   &      &  \\
}$$

And the summands of $T$ are:

$\xymatrix @! @=0.15pc{
      &   && & P_1 \\
    &  & P_2 \ar[rr]^{\beta}& & P_1 \\
    P_3 \ar[rr]^{\beta}& & P_2 & &   \\
    P_4\ar[rr]^{\beta^2}& & P_2& &   \\
     & &   && P_5. \\
}$
\end{ex}

Let us compute the tilting complex $T$ from the lemma. Denote by
$\psi _{\Gamma}: \{1,2,\dots,n\}  \backslash \phi_{\Gamma}^{-1}(0)
\rightarrow \{1,\dots,n-1\}$ the function which assigns to a label $i$
the label of the edge from a higher level which shares a common
vertex with the edge $i$. Let us assume that in the complexes $0
\rightarrow P_i \rightarrow 0$, $0 \rightarrow P_i \rightarrow P_j
\rightarrow 0,$ $P_i$ is concentrated in degree $0$, (note that we are
using the cohomological notation).

\begin{lem} Let $\Gamma$ be a Brauer tree of type $(n,t)$ labelled as described
above. Let $T$ be the complex from Lemma \ref{l5}, then $T$ is of the form
$T= \bigoplus_{i=1}^{n}T_{i}$, where

$$T_i=\left\{
            \begin{array}{ll}
              P_i, & \phi_{\Gamma}(i)=0,  \\
(P_i \xrightarrow{\beta^{i-\psi(i)}} P_{\psi(i)})[\phi_{\Gamma}(i)],
& \phi_{\Gamma}(i)\neq 0.
\end{array}
          \right.
$$

\end{lem}

\textbf{Proof.} Recall that the arrows in the Brauer star algebra
were denoted by $\beta$. Since $T= (\mu^{+}_n)^{\phi_{\Gamma}(n)}
\circ(\mu^{+}_{n-1})^{\phi_{\Gamma}(n-1)}\circ\dots\circ
(\mu^{+}_1)^{\phi_{\Gamma}(1)}(A)$, first we mutate the summand with
label $1$, then with label $2$ and so on, hence we can compute $T$
by induction on the label. $\phi_{\Gamma}(1)=0$, hence $T_1=P_1.$ If
the edge with label 1 is not incident to a leaf, then 2 is on the
edge incident to its non-exceptional end which is the previous edge
coming before the edge with label 1 in the cyclic ordering,
$\phi_{\Gamma}(2)=1$ and we should apply $\mu^{+}_2$ to $A$. Hence
$T_2= P_2 \overset{\beta}{\rightarrow} P_{1}$ concentrated in
degrees $-1$ and $0$, as desired. If the edge with label 1 is
incident to a leaf, then 2 is on the edge incident to the
exceptional vertex, $\phi_{\Gamma}(2)=0$, hence $T_2=P_2.$

Assume that $T_i$ is computed. If the edge with label $i+1$ is
incident to the exceptional vertex, then $\phi_{\Gamma}(i+1)=0$ and
$T_{i+1}=P_{i+1}$.

Assume the label $i+1$ is assigned to some edge not incident to the
exceptional vertex. Denote by
$x_1,x_2,\dots, x_{\phi_{\Gamma}(i+1)}=\psi(i+1)$ the edges from the
shortest path from the exceptional vertex to $i+1$ indexed by
numbers from $1$ to $\phi_{\Gamma}(i+1)$. By assumption
$$T^{r}= (\mu^{+}_i)^{\phi_{\Gamma}(i)}
\circ(\mu^{+}_{i-1})^{\phi_{\Gamma}(i-1)}\circ\dots\circ
(\mu^{+}_1)^{\phi_{\Gamma}(1)}(A)= \bigoplus_{k=1}^{i}T_{k} \oplus
\bigoplus_{k=i+1}^{n}P_{k}.$$
$$T_{x_1}=P_{x_{1}}, \mbox{ } T_{x_{k}} \simeq (P_{x_{k}}
\xrightarrow{\beta^{x_{k}-x_{k-1}}}
P_{x_{k-1}})[\phi_{\Gamma}(x_{k})].$$ We want to apply
$(\mu^{+}_{i+1})^{\phi_{\Gamma}(i+1)}$ to $T^{r}$, or equivalently
mutate $P_{i+1}$ along $T_{x_1}$,
$T_{x_2},\dots,T_{x_{\phi_{\Gamma}(i+1)}}$. Clearly $T^{r+1}_{i+1}
\simeq P_{i+1} \xrightarrow{\beta^{i+1-x_{1}}} P_{x_{1}}$. The
minimal left approximation of $T^{r+1}_{i+1}$ with respect to other
summands of $T^{r+1}$ is
$$\xymatrix{P_{i+1} \ar[d]^{\beta^{i+1-x_{2}}} \ar[r]^{\beta^{i+1-x_{1}}}
&P_{x_{1}} \ar[d]^{Id} \\
P_{x_{2}} \ar[r]^{\beta^{x_{2}-x_{1}}} & P_{x_{1}}.}$$ By Lemma 4
$T^{r+2}_{i+1} \simeq (P_{i+1} \xrightarrow{\beta^{i+1-x_{2}}}
P_{x_{2}})[2]$, and by iterated application of Lemma 4 we get the
desired result.\hfill\(\Box\)

\begin{rem}
Complexes similar to $T$ from Lemma \ref{l5} were already studied in the works of Schaps and Zakay-Illouz (see
for example \cite{SI2}).
\end{rem}

\subsection{Main result}

Let $A$ be a Brauer star algebra of type (n,t). Denote by $\mathcal{R}$
the subgroup of ${\rm TrPic}(A)$ generated by shift, ${\rm
Pic}(A)$ and equivalences $H_i$ (see formula (1)).

\begin{rem}

\emph{(a)} The subgroup $\mathcal{R}$ coincides with the subgroup
considered in \cite{SI}, it was also shown there that this subgroup
has an action of the braid group on the diagram $\widetilde{A}_{n-1}$,
the homomorphism is defined by sending half-twists to $H_i$'s. In
\cite{IM} this action was shown to be faithful for $t=1$.

\emph{(b)} There is an outer automorphism in ${\rm Pic}(A)$
corresponding to the rotation of the Brauer star and sending $H_i$
to $H_{i+1}$ by conjugation, so one can define $\mathcal{R}$ as a
subgroup of ${\rm TrPic}(A)$ generated by shift, ${\rm Pic}(A)$ and
equivalence $H_1$.

\emph{(c)} $H_i(A) \simeq (\mu^+_i)^2(A)$: for a complex $$P_i
\xrightarrow{\beta} P_{i-1} \xrightarrow{soc} P_{i-1}$$ there are a
triangle $$P_i \xrightarrow{f'} M' \rightarrow Cone(f') \rightarrow,
$$ where $M'=P_{i-1}$, $f'=\beta$ is a minimal left approximation of
$P_i$ with respect to ${\rm add}(\bigoplus_{j=1, j \neq i}^{n}P_j)$
and $Cone(f')\simeq P_i \xrightarrow{\beta} P_{i-1};$ and a triangle
$$Cone(f') \xrightarrow{f''} M'' \rightarrow Cone(f'') \rightarrow, $$
where $M''=P_{i-1}$, $f''=(0,soc)$ is a minimal left approximation
of $Cone(f')$ with respect to ${\rm add}(\bigoplus_{j=1, j \neq
i}^{n}P_j)$ and $$Cone(f'') \simeq (P_i \xrightarrow{\beta} P_{i-1}
\xrightarrow{soc} P_{i-1}).$$
\end{rem}

\begin{thm}
If $t>1$, then ${\rm TrPic}(A) = \mathcal{R}$.
\end{thm}

\textbf{Proof.} We only need to show that the
embedding of $\mathcal{R}$ into ${\rm  TrPic}(A)$ is surjective.

\emph{1st step: It is sufficient to prove that for any tilting
complex $T$ there is an element from $\mathcal{R}$ which sends
indecomposable projective $A$-modules to the summands of $T$.}

Indeed, any element $X$ from ${\rm TrPic}(A)$ restricts to some tilting complex
$T$ such that ${\rm End}_{D^b(A)}(T)^{op} \simeq A$; any other element
from ${\rm TrPic}(A)$ which restricts to the same complex, differs
from $X$ by an element from ${\rm Pic}(A)$, hence by an element from
$\mathcal{R}$.

\emph{2nd step: If $T=(\mu^{\pm}_{j_q})^2
 \circ\dots\circ (\mu^{\pm}_{j_{1}})^2(A)$, then there exists an element from $\mathcal{R}$ which sends
the indecomposable projective $A$-modules to the summands of $T$.}

Indeed, let $F$ be some autoequivalence of $D^b(A)$, let us compute
$F(H_i(P_j))$. Since $H_i(P_j)= P_j$ for $j \neq i, i-1$, we have
$F(H_i(P_j))=F(P_j)$ for $j \neq i, i-1$. $H_i(P_i) = P_{i-1}$,
hence $F(H_i(P_i))=F(P_{i-1})$. For the complex $H_i(P_{i-1})$ there
are triangles $$P_i \xrightarrow{f'} M' \rightarrow Cone(f')
\rightarrow $$ and $$Cone(f') \xrightarrow{f''} M'' \rightarrow
Cone(f'') \rightarrow, $$ where $f'$ and $f''$ are minimal left
approximations of $P_i$ and $Cone(f')$ with respect to ${\rm
add}(\bigoplus_{j=1, j \neq i}^{n}P_j)$, hence by Lemma $1$ for the
complex $F(H_i(P_{i-1}))$ there are triangles $$F(P_i)
\xrightarrow{g'} N' \rightarrow Cone(g') \rightarrow $$ and $$Cone(g')
\xrightarrow{g''} N'' \rightarrow Cone(g'') \rightarrow, $$ where
$g'=F(f')$ and $g''=F(f'')$ are minimal left approximations of
$F(P_i)$ and $Cone(g') \simeq F(Cone(f'))$ with respect to ${\rm
add}(\bigoplus_{j=1, j \neq i}^{n}F(P_j))$ and $F(H_i(P_{i-1}))
\simeq Cone(g'')$. So $F(H_i(P_{i-1}))$ is a double mutation of
$F(P_i)$ with respect to other summands of $F(A)$. Consider an
autoequivalence such that
$$H'_i(P_j)=\left\{
            \begin{array}{ll}
              P_j\mbox{ } \rightarrow \mbox{ } 0\mbox{ } \mbox{ }\rightarrow \mbox{ } 0, & j \neq i, i-1  \\
              P_{i}\mbox{ }\rightarrow \mbox{ } 0\mbox{ } \mbox{ }\rightarrow \mbox{ } 0, & j = i-1 \\
              P_i \xrightarrow{soc} P_{i} \xrightarrow{\beta} P_{i-1}, & j
              =i.
            \end{array}
          \right.$$
Note that $H'_i(A) \simeq (\mu^-_{i-1})^2(A)$, i.e. for $H'_i(P_i)$
there are two triangles from the definition of the right mutation.
Analogously to the previous argument we get that to apply some
autoequivalence $G$ to $H'_i(A)$ is the same as to compute a double
right mutation of $G(P_{i-1})$ with respect to other summands of
$G(A)$. So $H'_i(H_i(P_j))=P_j$ and $H_i(H'_i(P_j))=P_j$ for any
$j$, hence the action of $(H_i)^{-1}$ on the projective $A$-modules
coincides with the action of $H'_i$.

Assume that some tilting complex $T$ can be obtained from $A$
applying squares of mutations, i.e. $T=(\mu^{\pm}_{j_q})^2
 \circ\dots\circ (\mu^{\pm}_{j_{1}})^2(A)$. By induction on
$q$ we get that there is an element from $\mathcal{R}$, which sends
indecomposable projective $A$-modules to the summands of $T$.

\emph{3rd step:} Assume that $T$ is concentrated in non-positive degrees. By results
of Aihara~\cite{Ai2} $T=\mu^{+}_{i_s}
 \circ\dots\circ \mu^{+}_{i_{r+1}} \circ  \mu^{+}_{i_r} \circ
 \mu^{+}_{i_{r-1}} \circ\dots \circ \mu^{+}_{i_2} \circ
 \mu^{+}_{i_1}(A)$ for some $(i_1, i_2, \dots, i_s)$.
Denote by $T^r:=\mu^{+}_{i_r} \circ \mu^{+}_{i_{r-1}} \circ\dots \circ
\mu^{+}_{i_2} \circ \mu^{+}_{i_1}(A)$ and by $\Gamma^r$ the Brauer
tree of ${\rm End}_{D^b(A)}(T^r)^{op}$. Note that if the Brauer star is
labelled in the standard way, $\Gamma^r$ has a natural labelling of
edges associated to the series of mutations $\mu^{+}_{i_r} \circ
\mu^{+}_{i_{r-1}} \circ\dots \circ \mu^{+}_{i_2} \circ
\mu^{+}_{i_1}(A)$, which may not coincide with the standard
labelling.

By Lemma 5 if $A$ is a Brauer star algebra with the standard
labelling, then the opposite ring of the endomorphism ring of
$(\mu^{+}_n)^{\phi_{\Gamma^r}(n)}
\circ(\mu^{+}_{n-1})^{\phi_{\Gamma^r}(n-1)}\circ\dots\circ
(\mu^{+}_1)^{\phi_{\Gamma^r}(1)}(A)$ is a Brauer tree algebra
associated to $\Gamma^r$ with the standard labelling. If the natural
labelling $\rho_r$ of tree $\Gamma^r$ is not standard, then there is
some permutation $\tau_r$ one needs to apply to the standard
labelling of $\Gamma^r$ to obtain $\rho_r$. Applying $\tau_r$ to the
standard labelling of the Brauer star we obtain the Brauer star with
some labelling, which will also be denoted by  $\rho_r$. Applying
$\tau_r$ to the indices of $(\mu^{+}_n)^{\phi_{\Gamma^r}(n)}
\circ(\mu^{+}_{n-1})^{\phi_{\Gamma^r}(n-1)}\circ\dots\circ
(\mu^{+}_1)^{\phi_{\Gamma^r}(1)}(A)$ one gets a series of mutations
needed to obtain the Brauer tree ${\Gamma^r}$ with the labelling
$\rho_r$ from the Brauer star with the labelling $\rho_r$. We will
denote this series of mutations by $\tilde{\mu}^r$ for the natural
labelling $\rho_r$ of ${\Gamma^r}$. By
$\tilde{\mu}^{-r}:=(\tilde{\mu}^{r})^{-1}$ we will denote the series
of mutations
$(\mu^{-}_{\tau_r(1)})^{\phi_{\Gamma^r}(\tau_r(1))}\circ\dots\circ(\mu^{-}_{\tau_r(n-1)})^{\phi_{\Gamma^r}(\tau_r(n-1))}\circ(\mu^{-}_{\tau_r(n)})^{\phi_{\Gamma^r}(\tau_r(n))}$.
If $A^{\rho_r}$ is a Brauer star algebra with the labelling $\rho_r$
and $B^r$ is the Brauer tree algebra of ${\Gamma^r}$ with natural
labelling, then $\tilde{\mu}^{-r} \circ \tilde{\mu}^{r}(A^{\rho_r})
= A^{\rho_r}$ and $\tilde{\mu}^{r} \circ \tilde{\mu}^{-r}(B^r) =
B^r$.
$$T=\mu^{+}_{i_s}
 \circ\dots\circ \mu^{+}_{i_{r+1}} \circ  \mu^{+}_{i_r} \circ
 \mu^{+}_{i_{r-1}} \circ\dots \circ \mu^{+}_{i_2} \circ \mu^{+}_{i_1}(A)=$$
$$=\mu^{+}_{i_s}
 \circ\dots\circ \tilde{\mu}^{r+1} \circ
 \tilde{\mu}^{-(r+1)}\circ \mu^{+}_{i_{r+1}}
 \circ \tilde{\mu}^{r} \circ \tilde{\mu}^{-r}\circ
 \mu^{+}_{i_r}\circ \tilde{\mu}^{r-1} \circ \tilde{\mu}^{-(r-1)}
 \circ  \mu^{+}_{i_{r-1}} \circ\dots $$
$$\circ \tilde{\mu}^{2} \circ \tilde{\mu}^{-2} \circ
\mu^{+}_{i_2} \circ \tilde{\mu}^{1} \circ \tilde{\mu}^{-1} \circ
\mu^{+}_{i_1}(A).$$ By the 2nd step of the proof of the theorem
follows from the fact that any tilting complex $T$ can be obtained
from $A$ applying squares of mutations and the shift, thus we have
proved the following:

\emph{It is sufficient to express as series of squares of mutations
of $A$ only complexes of the form $\tilde{\mu}^{-(r+1)}\circ
\mu^{+}_{i_{r+1}}\circ \tilde{\mu}^{r}(A^{\rho_r})$,
$\tilde{\mu}^{-1} \circ \mu^{+}_{i_1}(A)$ and $\mu^{+}_{i_{s}}\circ
\tilde{\mu}^{s-1}(A^{\rho_{s-1}})$.}

Note that, the opposite rings of the endomorphism rings of these
tilting complexes are isomorphic to the Brauer star algebra.

\emph{4th step: The complexes of the form $\tilde{\mu}^{-(r+1)}\circ
\mu^{+}_{i_{r+1}}\circ \tilde{\mu}^{r}(A^{\rho_r})$,
$\tilde{\mu}^{-1} \circ \mu^{+}_{i_1}(A)$ and $\mu^{+}_{i_{s}}\circ
\tilde{\mu}^{s-1}(A^{\rho_{s-1}})$ can be expressed as series of squares of mutations
of $A$.}

Indeed, $\tilde{\mu}^{-1} \circ \mu^{+}_{i_1}(A)=A$ and
$\mu^{+}_{i_{s}}\circ \tilde{\mu}^{s-1}(A^{\rho_{s-1}}) =
(\mu^{+}_{i_{s}})^2(A^{\rho_{s-1}})$, here we use the fact that
$t>1$ and that the exceptional vertex is distinguished by the
multiplicity. So it remains to check the condition for
$\tilde{\mu}^{-(r+1)}\circ \mu^{+}_{i_{r+1}}\circ
\tilde{\mu}^{r}(A^{\rho_r})$ for any tree ${\Gamma^r}$ and any
mutation $\mu^{+}_{i_{r+1}}$. Applying $\tau_r^{-1}$we are going to
perform a suitable series of mutations with the standard Brauer
star. This is done in sections \ref{1} and \ref{2}.

Thus the theorem follows.\hfill\(\Box\)

\begin{rem}
The question of the relations in $\mathcal{R}$ remains open.
\end{rem}

\subsection{Mutation of type \MakeUppercase{\romannumeral1}}\label{1}

Note that these computations are also valid for $t=1$. Often we will
omit the degrees is which the complexes are concentrated, but will
take them into account.

Without loss of generality assume that ${\Gamma^r}$ is a tree with a
standard labelling, so $\tilde{\mu}^{r}(A)$ is a series of mutations
described in section 3.1. Assume also that $j:=i_{r+1}$ is an edge
incident to a vertex $x$ of degree one in ${\Gamma^r}$, so
$\mu^{+}_{j}$ is a mutation of type \MakeUppercase{\romannumeral1}.
In order to compute $\tilde{\mu}^{-(r+1)}\circ  \mu^{+}_{j}\circ
\tilde{\mu}^{r}(A)$ we are going to consider the following three
cases: 1) $x$ is not the root and there is an edge of the same level
as $j$, following $j$ in the cyclic ordering; 2) $x$ is not the root
and there is no edge of the same level as $j$, following $j$ in the
cyclic ordering; 3) $x$ is the root.

\textbf{Case 1:} $x$ is not the root and there is an edge $l$ of the same level as
$j$, following $j$ in the cyclic ordering in ${\Gamma^r}$.

\emph{Claim:} \emph{${\Gamma^{r+1}}$ is a tree with the standard
labelling, $\mu^{+}_{j}\circ \tilde{\mu}^{r}(A)$ is isomorphic to
the standard tilting complex from Lemma 6, associated to this tree,
and hence $\tilde{\mu}^{-(r+1)}\circ  \mu^{+}_{j}\circ
\tilde{\mu}^{r}(A) \simeq A$.}

\emph{Proof.} It is clear that ${\Gamma^{r+1}}$ is a tree with the
standard labelling. Assume that $j$ is not incident to the root,
denote by $q$ the edge of the higher level incident to both $j$ and
$l$ in ${\Gamma^r}$. Since $\phi_{\Gamma^r}(j) = \phi_{\Gamma^r}(l)$
and $\psi(j)=\psi(l)=q$ in ${\Gamma^r}$ the minimal left
approximation $f$ of $T^r_j$ has the form:

$$\xymatrix{T^r_j \simeq \ar[d]^{f} & P_{j} \ar[d]^{\beta^{j-l}} \ar[r]^{\beta^{j-q}} &P_{q} \ar[d]^{{\rm Id}} \\
T^r_l \simeq & P_{l} \ar[r]^{\beta^{l-q}} & P_{q},}$$ since ${\rm
dim_k{\rm Hom}}_{D^b(A)}(T^r_j,T^r_l)=1$ and any non-zero morphism
is a minimal left approximation. $T^{r+1}_{j}\simeq(P_{j}
\xrightarrow{\beta^{j-l}} P_{l})[\phi_{\Gamma^r}(j)+1]$   by Lemma
4, as desired, since $\phi_{\Gamma^{r+1}}(j) = \phi_{\Gamma^r}(j)+1$
and $\psi(j)=l$ in ${\Gamma^{r+1}}$.

If $j$ is incident to the root in ${\Gamma^{r}}$, $l$ is also
incident to the root in ${\Gamma^{r}}$ and $T^{r+1}_{j}=(P_{j}
\overset{\beta^{j-l}}{\rightarrow} P_{l})[1]$.\hfill\(\Box\)

\textbf{Case 2:} $x$ is not the root and there is no edge of the same level as
$j$, following $j$ in the cyclic ordering (in ${\Gamma^r}$). In
order to simplify the notation set $\phi_{\Gamma}=\phi_{\Gamma^r}.$

\emph{Claim:} \emph{ $\tilde{\mu}^{-(r+1)}\circ \mu^{+}_{j}\circ
\tilde{\mu}^{r}(A) \simeq (\mu^+_j)^2(A)$.}

\emph{Proof.} Let $T^r = \bigoplus_{i=1}^{n}T_i$ be the tilting
complex from Lemma 6, and let $j$ be not incident to the root.
Denote by $x_1,x_2,\dots,x_{\phi_{\Gamma}(j)}=\psi(j)=j-1$ the edges
which belong to the shortest path from the exceptional vertex to the
edge labelled $j$, indexed by numbers from $1$ to
$\phi_{\Gamma}(j)$. $T^{r+1}_j$ is the cone of the following
morphism (since ${\rm dim_k{\rm Hom}}_{D^b(A)}(T^r_j,T^r_{j-1})=1$
any non-zero morphism is a minimal left approximation)
$$\xymatrix   {
 P_j \ar[d] \ar[r]^{\beta} & P_{j-1}\ar[d]^{soc} \ar[r] & 0 \ar[d] \\
   0 \ar[r] & P_{j-1}\ar[r]^{\beta^{a}} & P_{x_{\phi_{\Gamma}(j)-1}}, \\
}
$$
where $soc$ is the morphism with the image equal to the socle of
$P_{j-1}$, and $a=j-1-x_{\phi_{\Gamma}(j)-1}$. Hence $$T^{r+1}_j
\simeq P_j \xrightarrow{\beta}P_{j-1} \xrightarrow{soc} P_{j-1}
\xrightarrow{\beta^{a}} P_{x_{\phi_{\Gamma}(j)-1}},$$ where
$P_{x_{\phi_{\Gamma}(j)-1}}$ is concentrated in the same degree as
in $T^r$. If the edge $j-1$ is incident to the root, then
$P_{x_{\phi_{\Gamma}(j)-1}}=0$ and $$T^{r+1}_j \simeq P_j
\xrightarrow{\beta}P_{j-1} \xrightarrow{soc} P_{j-1}.$$ Let us apply
the part of the series of mutations $\tilde{\mu}^{-(r+1)}$, which
contains the mutation with indices greater than $j$ and with index
$j-1$. For the corresponding $r'$ we get that $$T^{r'} \simeq
\bigoplus_{i=1}^{j-2} T_i \oplus (P_j \xrightarrow{\beta}P_{j-1}
\xrightarrow{soc} P_{j-1} \xrightarrow{\beta^a}
P_{x_{\phi_{\Gamma}(j)-1}}) \oplus P_{j-1} \bigoplus_{i=j+1}^{n}
P_i.$$ Let us apply $\mu^{-}_{j}$: $T^{r'+1}_j$ is the cone of the
following morphism shifted by $-1$,
$$\xymatrix  {
 0\ar[r] \ar[d] & 0\ar[r] \ar[d] & 0\ar[r] \ar[d] &P_{x_{\phi_{\Gamma}(j)-1}}\ar[d]^{{\rm Id}} \ar[r]^{\beta^{b}} &P_{x_{\phi_{\Gamma}(j)-2}} \ar[d] \\
  P_j  \ar[r]^{\beta} & P_{j-1}\ar[r]^{soc} & P_{j-1}\ar[r]^{\beta^{a}} & P_{x_{\phi_{\Gamma}(j)-1}}\ar[r]& 0,\\
}
$$ where $b=x_{\phi_{\Gamma}(j)-1}-x_{\phi_{\Gamma}(j)-2}$.
By Lemma 4 $$T^{r'+1}_j \simeq P_j \xrightarrow{\beta}P_{j-1}
\xrightarrow{soc} P_{j-1} \xrightarrow{\beta^{a+b}}
P_{x_{\phi_{\Gamma}(j)-2}},$$ where $P_{x_{\phi_{\Gamma}(j)-2}}$ is
concentrated in the same degree as in $T^r$. By iterated application
of Lemma 4 $$T_j \simeq P_j \xrightarrow{\beta}P_{j-1}
\xrightarrow{soc} P_{j-1},$$ where $P_{j}$ is concentrated in degree $-2$
(from here on by $T$ we will denote the complex
$\tilde{\mu}^{-(r+1)}\circ \mu^{+}_{j}\circ \tilde{\mu}^{r}(A)$ and
by $T_j$ we will denote its summand corresponding to the edge with
label $j$).

Applying the rest of the series of mutations $\tilde{\mu}^{-(r+1)}$,
we get $T_i \simeq P_i$ for $i<j-1$, hence
$(\mu^+_j)^2(A)$.\hfill\(\Box\)

\textbf{Case 3:} $x$ is the root, $j=1$.

Let $l$ be the edge following $j$ in the cyclic ordering around the
non-root vertex incident to $j$ in $\Gamma^r$.

$\Gamma={\Gamma^{r}}$ and ${\Gamma^{r+1}}$ are roughly of the
following form (labels $\{2,\dots,l-1\}$ and $\{l+1,\dots,n\}$ are on
arbitrary trees):

$$
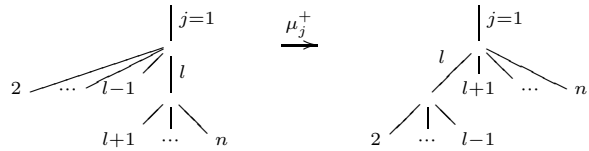
\xymatrix @! @=0.01pc {
  & & & \ar@{-}[d]^{j=1}  & & & & && \ar@{-}[d]^{j=1}   &&\\
    & & & \ar@{-}[d]^{l}  \ar@{-}[dlll] \ar@{-}[dll] \ar@{-}[dl]   && \ar[r]^{\mu^{+}_{j}}&&  & & \ar@{-}[d] \ar@{-}[dl]_{l} \ar@{-}[dr] \ar@{-}[d] \ar@{-}[drr]  &&\\
     ^{2} & ^{\cdots} & ^{l-1} & \ar@{-}[dl]   \ar@{-}[dr] \ar@{-}[d] & && & & \ar@{-}[d] \ar@{-}[dl] \ar@{-}[dr]& ^{l+1} & ^{\cdots}&^{n}\\
   &    & ^{l+1} &  ^{\cdots} & ^{n}&&&^{2}&^{\cdots}&   ^{l-1} & & \\
}$$

\emph{Claim:} \emph{ $\tilde{\mu}^{-(r+1)}\circ \mu^{+}_{j}\circ
\tilde{\mu}^{r}(A) \simeq ((\mu^-_1)^2 \circ (\mu^-_2)^2 \circ \dots
\circ (\mu^-_{l-1})^2(A))[1]$.}

\emph{Proof.} Let $T^r$ be the complex from Lemma 6, corresponding
to the tree $\Gamma^{r}$. Thus $T^r_1=P_1$ and $T^r_l \simeq
P_l\xrightarrow{\beta^{l-1}} P_1$. The minimal left approximation of
$P_1$ is the morphism

$$\xymatrix  { 0 \ar[r] \ar[d] & P_1 \ar[d]^{{\rm Id}}\\
P_l \ar[r]^{\beta^{l-1}} & P_1.\\
}
$$
Its cone is $P_l[1]$.

The levels of edges in $\Gamma^{r+1}$ are: $\phi_{\Gamma^{r+1}}(2) =
\phi_{\Gamma}(2)+1,\dots, \phi_{\Gamma^{r+1}}(l-1) =
\phi_{\Gamma}(l-1)+1$, $\phi_{\Gamma^{r+1}}(l+1) =
\phi_{\Gamma}(l+1)-1,\dots, \phi_{\Gamma^{r+1}}(n) =
\phi_{\Gamma}(n)-1$, hence the sequence of mutations
$\tilde{\mu}^{-(r+1)}$ is of the form $$\mu^-_l \circ
(\mu^-_{2})^{\phi_{\Gamma}(2)+1}\circ \dots \circ
(\mu^-_{l-1})^{\phi_{\Gamma}(l-1)+1} \circ
(\mu^-_{l+1})^{\phi_{\Gamma}(l+1)-1} \circ\dots \circ
(\mu^-_{n})^{\phi_{\Gamma}(n)-1}.$$ Let us apply this sequence to the
complex $T^{r+1}$.

Application of $(\mu^-_{n})^{\phi_{\Gamma}(n)-1}$: note that the
edge $n$ is pendant; if $n$ is not incident to $j$, then after the
application of $(\mu^-_{n})^{\phi_{\Gamma}(n)-2}$ the edge $n$
becomes incident to $j$, the complex corresponding to that edge is
$P_n \xrightarrow{\beta^{n-l}} P_l$. The minimal right approximation
of this complex is the morphism from $P_l[1]$ induced by $P_l
\xrightarrow{{\rm Id}} P_l$, by Lemma 4 its cone is isomorphic to
$P_n[2]$, hence $T_n \simeq P_n[1]$.

Similarly applying $(\mu^-_{l+1})^{\phi_{\Gamma}(l+1)-1} \circ\dots
\circ (\mu^-_{n-1})^{\phi_{\Gamma}(n-1)-1}$ we get that $$T_{n-1}
\simeq P_{n-1}[1],\dots, T_{l+1} \simeq P_{l+1}[1].$$

Application of $(\mu^-_{l-1})^{\phi_{\Gamma}(l-1)+1} $: note that
the edge $l-1$ is pendant; if $l-1$ is not incident to $l$, then
after the application of $(\mu^-_{l-1})^{\phi_{\Gamma}(l-1)-1}$, the
edge $l-1$ becomes incident to $l$, the complex corresponding to
that edge is $P_{l-1} \xrightarrow{\beta^{l-2}} P_j$. The minimal
right approximation of this complex is the morphism $(\beta,{\rm
Id})$ from $P_{l} \xrightarrow{\beta^{l-1}} P_j$. By Lemma 4 its
cone is isomorphic to $P_{l} \xrightarrow{\beta} P_{l-1}$. The
minimal right approximation of the shift of this complex is the
morphism from $P_l[1]$ induced by $P_l \xrightarrow{soc} P_l$, its
cone is isomorphic to $$P_l \xrightarrow{soc}P_{l}
\xrightarrow{\beta} P_{l-1},$$ after the shift $P_{l-1}$ is
concentrated in $1$, hence $$T_{l-1} \simeq (P_l
\xrightarrow{soc}P_{l} \xrightarrow{\beta} P_{l-1})[1].$$

Similarly applying $(\mu^-_{2})^{\phi_{\Gamma}(2)+1}\circ \dots \circ
(\mu^-_{l-2})^{\phi_{\Gamma}(l-2)+1}$ we get that
$$T_{l-2} \simeq  (P_l \xrightarrow{soc}P_{l} \xrightarrow{\beta^2} P_{l-2})[1],$$
$$\dots$$
$$T_{2} \simeq  (P_l \xrightarrow{soc}P_{l} \xrightarrow{\beta^{l-2}}
P_{2})[1].$$

Let us apply $\mu^-_l$: the minimal right approximation of $P_l
\xrightarrow{\beta^{l-j}} P_j$ is the morphism from $P_l[1]$,
induced by $P_l \xrightarrow{soc} P_l$, its cone is isomorphic to
$$P_l \xrightarrow{soc}P_{l} \xrightarrow{\beta^{l-j}} P_{j},$$ after
the shift $P_{j}$ is concentrated in $1$, hence $$T_{l} \simeq (P_l
\xrightarrow{soc}P_{l} \xrightarrow{\beta^{l-1}} P_{j})[1] = (P_l
\xrightarrow{soc}P_{l} \xrightarrow{\beta^{l-1}} P_{1})[1].$$

Thus $T \simeq ((\mu^-_1)^2 \circ (\mu^-_2)^2 \circ \dots \circ
(\mu^-_{l-1})^2(A))[1]$.\hfill\(\Box\)

Note that the indices $j$ and $l$ have switched places.

\subsection{Mutation of type \MakeUppercase{\romannumeral2}}\label{2}

Assume that both ends of $j:=i_{r+1}$ are vertices of degree $>1$ in
${\Gamma^r}$, so $\mu^{+}_{j}$ is a mutation of type
\MakeUppercase{\romannumeral2}. In order to compute
$\tilde{\mu}^{-(r+1)}\circ \mu^{+}_{j}\circ \tilde{\mu}^{r}(A)$ we
are going to consider the following three cases: 1) $j$ is not
incident to the root and there is an edge of the same level as $j$,
following $j$ in the cyclic ordering; 2) $j$ is incident to the
root; 3) $j$ is not incident to the root and there is no edge of the
same level as $j$, following $j$ in the cyclic ordering.

\textbf{Case 1:} $j$ is not incident to the root and there is an edge $l$ of the
same level as $j$, following $j$ in the cyclic ordering. The edge of
the next level, following $j$ in the cyclic ordering is denoted by
$m$.

$\Gamma={\Gamma^{r}}$ and ${\Gamma^{r+1}}$ are roughly of the
following form (labels $\{j+1,\dots,m-1\}$ and $\{m+1,\dots,h\}$ are on
arbitrary trees, $l$ is not necessarily a pendant edge):

$$
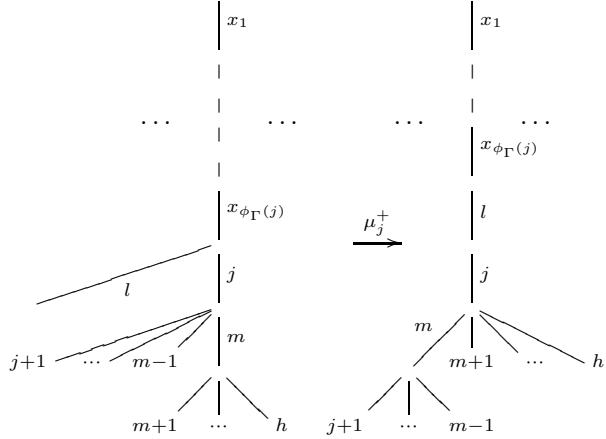
\xymatrix @! @=0.01pc {
  &  &  & \ar@{-}[d]^{x_1} &   &  &  & \ar@{-}[d]^{x_1} & & \\
  &&& \ar@{-}@{--}[d]  &   &&& \ar@{-}@{--}[d]  & & \\
&&\dots& \ar@{-}@{--}[d]  & \dots   &&\dots& \ar@{-}[d]^{x_{\phi_{\Gamma}(j)}}  &\dots  & \\
  & & & \ar@{-}[d]^{x_{\phi_{\Gamma}(j)}}  &  & & & \ar@{-}[d]^{l}  &&\\
  & & & \ar@{-}[d]^{j} \ar@{-}[dlll]^{l} & & \ar[r]^{\mu^{+}_{j}} & & \ar@{-}[d]^{j}  &&\\
    & & & \ar@{-}[d]^{m}  \ar@{-}[dlll] \ar@{-}[dll] \ar@{-}[dl]   & & & & \ar@{-}[d] \ar@{-}[dl]_{m} \ar@{-}[dr] \ar@{-}[d] \ar@{-}[drr]  &&\\
     ^{j+1} & ^{\cdots} & ^{m-1} & \ar@{-}[dl]   \ar@{-}[dr] \ar@{-}[d] & & &  \ar@{-}[d] \ar@{-}[dl] \ar@{-}[dr]& ^{m+1} & ^{\cdots}&^{h}\\
   &    & ^{m+1} &  ^{\cdots} & ^{h}&^{j+1}&^{\cdots}&   ^{m-1} & & \\
}$$

As before let us denote by
$x_1,x_2,\dots,x_{\phi_{\Gamma}(j)}=\psi(j)$ the edges which belong to
the shortest path from the exceptional vertex to the edge labelled
$j$, indexed by numbers from $1$ to $\phi_{\Gamma}(j)$.

\emph{Claim:} \emph{ $\tilde{\mu}^{-(r+1)}\circ \mu^{+}_{j}\circ
\tilde{\mu}^{r}(A) \simeq (\mu^-_{j})^2\circ (\mu^-_{j+1})^2\circ
\dots \circ (\mu^-_{m-1})^2(A)$.}

\emph{Proof.} Let $T^r$ be the complex from Lemma 6 associated to
the tree $\Gamma^{r}$. Since $\phi_{\Gamma^r}(j) =
\phi_{\Gamma^r}(l)$, $\phi_{\Gamma^r}(m) = \phi_{\Gamma^r}(j)+1$ and
$\psi(j)=\psi(l)=x_{\phi_{\Gamma}(j)}$ in ${\Gamma^r}$, the minimal
left approximation $f$ of $T^r_j$ is of the form:

{\small $$\xymatrixcolsep{4pc}\xymatrix @!  {
 0 \ar[d] \ar[r]& P_{j}\ar[d]^{ \left(
                                    \begin{smallmatrix}
                                      {\rm Id} \\
                                      \beta^{j-l} \\
                                    \end{smallmatrix}
                                  \right)
 } \ar[r]^{\beta^{j-x_{\phi_{\Gamma}(j)}}} & P_{x_{\phi_{\Gamma}(j)}} \ar[d]^{{\rm Id}} \\
   P_m  \ar[r]^(.45){\small \left(
                                    \begin{smallmatrix}
                                      \beta^{m-j} \\
                                       0\\
                                    \end{smallmatrix}
                                  \right)
 } & P_{j} \oplus P_l \ar[r]^(.55){(0,\beta^{l-x_{\phi_{\Gamma}(j)}})} & P_{x_{\phi_{\Gamma}(j)}}, \\
}
$$}
since the edge corresponding to $T^r_j$ does not belong to the
exceptional cycle, hence a morphism consisting of non-zero morphisms
between the respective summands is a minimal left approximation.

The following diagram is commutative:

{\tiny $$ \xymatrixcolsep{2.5pc}\xymatrixrowsep{5pc}\xymatrix { &
P_m \ar[d]^{ \left(
                                    \begin{smallmatrix}
                                      {\rm Id} \\
                                      -\beta^{m-j} \\
                                    \end{smallmatrix}
                                  \right)
 }  \ar[rr]^{-\beta^{m-l}}&& P_{l}\ar[d]^{ \left(
                                    \begin{smallmatrix}
                                      0 \\
                                      {\rm Id}\\
                                      -\beta^{l-x_{\phi_{\Gamma}(j)}} \\
                                    \end{smallmatrix}
                                  \right)
 } \ar[rr] && 0 \ar[d] \\
Cone(f)= &   P_m \oplus P_j \ar[d]^{({\rm Id}, 0)}  \ar[rr]^(.45){
\left(
                                    \begin{smallmatrix}
                                      \beta^{m-j}& {\rm Id} \\
                                       0 & \beta^{j-l}\\
                                       0 & -\beta^{j-x_{\phi_{\Gamma}(j)}}\\
                                    \end{smallmatrix}
                                  \right)
 } && P_{j} \oplus P_l \oplus P_{x_{\phi_{\Gamma}(j)}}  \ar[d]^{(-\beta^{j-l}, {\rm Id}, 0)} \ar[rr]^(.6){(0,\beta^{l-x_{\phi_{\Gamma}(j)}}, {\rm Id})} && P_{x_{\phi_{\Gamma}(j)}} \ar[d] \\
 &  P_m \ar[rr]^{-\beta^{m-l}}  && P_{l}  \ar[rr] && 0, \\
}
$$}
hence $(P_m \xrightarrow{\beta^{m-l}}  P_{l}) \simeq (P_m
\xrightarrow{-\beta^{m-l}}  P_{l})$ is a summand of $Cone(f)$.
Clearly $Cone(f)$ is homotopic to $P_m \xrightarrow{\beta^{m-l}}
P_{l}$, since it is an indecomposable summand of a tilting complex.

The levels of edges in $\Gamma^{r+1}$ are: $\phi_{\Gamma^{r+1}}(j) =
\phi_{\Gamma}(j)+1$, $\phi_{\Gamma^{r+1}}(j+1) =
\phi_{\Gamma}(j+1)+2,\dots, \phi_{\Gamma^{r+1}}(m-1) =
\phi_{\Gamma}(m-1)+2$, $\phi_{\Gamma^{r+1}}(m) =
\phi_{\Gamma}(m)+1$, $\phi_{\Gamma^{r+1}}(m+1) =
\phi_{\Gamma}(m+1), \dots, \phi_{\Gamma^{r+1}}(h) =
\phi_{\Gamma}(h)$, hence $$\tilde{\mu}^{-(r+1)} =
(\mu^-_2)^{\phi_{\Gamma}(2)} \circ \dots \circ
(\mu^-_{j-1})^{\phi_{\Gamma}(j-1)} \circ
(\mu^-_{j})^{\phi_{\Gamma}(j)+1}
\circ(\mu^-_{m})^{\phi_{\Gamma}(m)+1}
\circ$$ $$\circ(\mu^-_{j+1})^{\phi_{\Gamma}(j+1)+2} \circ \dots \circ
(\mu^-_{m-1})^{\phi_{\Gamma}(m-1)+2} \circ
(\mu^-_{m+1})^{\phi_{\Gamma}(m+1)} \circ \dots \circ
(\mu^-_{h})^{\phi_{\Gamma}(h)}\circ $$ $$\circ
(\mu^-_{h+1})^{\phi_{\Gamma}(h+1)} \circ \dots \circ
(\mu^-_{n})^{\phi_{\Gamma}(n)}.$$ Let us apply this sequence to the
complex $T^{r+1}$.

Application of $(\mu^-_{h+1})^{\phi_{\Gamma}(h+1)} \circ \dots \circ
(\mu^-_{n})^{\phi_{\Gamma}(n)}$ gives $T_n \simeq P_n, \dots, T_{h+1}
\simeq P_{h+1}$.

Application of $(\mu^-_{h})^{\phi_{\Gamma}(h)}$: note that the edge
$h$ is pendant; if $h$ is not incident to $j$, then after the
application of $(\mu^-_{h})^{\phi_{\Gamma}(h)-\phi_{\Gamma}(m)-1}$
the edge $h$ becomes incident to $j$, the complex corresponding to
that edge is $P_h \xrightarrow{\beta^{h-m}} P_m$. The minimal right
approximation of this complex is the morphism from $P_m
\xrightarrow{\beta^{m-l}} P_{l}$ induced by $P_m \xrightarrow{{\rm
Id}} P_m$, by Lemma 4 its cone is isomorphic to $P_h
\xrightarrow{\beta^{h-l}} P_{l}$. The minimal right approximation of
the shift of $P_h \xrightarrow{\beta^{h-l}}  P_{l}$ is the morphism
from $P_l\xrightarrow{\beta^{l-x_{\phi_{\Gamma}(j)}}}~P_{x_{\phi_{\Gamma}(j)}}$ induced by $P_l \xrightarrow{{\rm Id}}
P_l$, by Lemma 4 its cone is isomorphic to $P_h
\xrightarrow{\beta^{h-l}} P_{x_{\phi_{\Gamma}(j)}}$, by iterated
application of Lemma 4 we get: $T_h \simeq P_h$.

Similarly, applying $(\mu^-_{m+1})^{\phi_{\Gamma}(m+1)} \circ \dots
\circ (\mu^-_{h-1})^{\phi_{\Gamma}(h-1)}$ we get: $$T_{h-1} \simeq
P_{h-1}, \dots, T_{m+1} \simeq P_{m+1}.$$

Application of $(\mu^-_{m-1})^{\phi_{\Gamma}(m-1)+2}$: note that the
edge $m-1$ is pendant; if $m-1$ is not incident to $m$, then after
the application of
$(\mu^-_{m-1})^{\phi_{\Gamma}(m-1)-\phi_{\Gamma}(m)}$ the edge $m-1$
becomes incident to $m$, the complex corresponding to that edge is
$P_{m-1} \xrightarrow{\beta^{m-1-j}} P_j$. The minimal right
approximation of this complex is the morphism $(\beta, {\rm Id})$
from $P_m \xrightarrow{\beta^{m-j}} P_{j}$, its cone is isomorphic
to $P_m~\xrightarrow{\beta}~P_{m-1}$. The minimal right
approximation of the shift of $P_m~\xrightarrow{\beta}~P_{m-1}$ is
the morphism $(soc, 0)$ from $P_m \xrightarrow{\beta^{m-l}} P_{l}$,
its cone is isomorphic to $$P_m
\xrightarrow{{\left( \begin{smallmatrix} -\beta^{m-l} \\  soc\\
\end{smallmatrix} \right)}}P_{l} \oplus P_{m} \xrightarrow{(0, \beta)}
P_{m-1}.$$ The minimal right approximation of this complex is the
morphism from $P_l \xrightarrow{\beta^{l-x_{\phi_{\Gamma}(j)}}}
P_{x_{\phi_{\Gamma}(j)}}$ induced by $P_l \xrightarrow{{\rm Id}}
P_l$, by a computation analogous to Lemma 4 we get that the cone is
isomorphic to $$P_m
\xrightarrow{{\left( \begin{smallmatrix} -\beta^{m-x_{\phi_{\Gamma}(j)}} \\  soc\\
\end{smallmatrix} \right)}}P_{x_{\phi_{\Gamma}(j)}} \oplus P_{m} \xrightarrow{(0, \beta)}
P_{m-1}.$$ Applying the remaining degrees of $\mu^-_{m-1}$ we get:
$$T_{m-1} \simeq P_m \xrightarrow{soc} P_{m} \xrightarrow{\beta} P_{m-1}.$$

Similarly applying $(\mu^-_{j+1})^{\phi_{\Gamma}(j+1)+2} \circ \dots
\circ (\mu^-_{m-2})^{\phi_{\Gamma}(m-2)+2}$ we get:
$$T_{m-2} \simeq P_m \xrightarrow{soc} P_{m} \xrightarrow{\beta^2} P_{m-2},$$
$$\dots$$
$$T_{j+1} \simeq P_m \xrightarrow{soc} P_{m} \xrightarrow{\beta^{m-j-1}} P_{j+1}.$$

Application of $(\mu^-_{m})^{\phi_{\Gamma}(m)+1}$: the minimal right
approximation of $P_m~\xrightarrow{\beta^{m-j}}~P_{j}$ is the
morphism $(soc, 0)$ from $P_m \xrightarrow{\beta^{m-l}} P_{l}$,
analogously to the previous case we get:
$$T_{m} \simeq P_m \xrightarrow{soc} P_{m}
\xrightarrow{\beta^{m-j}} P_{j}.$$

Application of $(\mu^-_{j})^{\phi_{\Gamma}(j)+1}$: the minimal right
approximation of $P_m \xrightarrow{\beta^{m-l}} P_{l}$ is the
morphism from $P_l \xrightarrow{\beta^{l-x_{\phi_{\Gamma}(j)}}}
P_{x_{\phi_{\Gamma}(j)}}$ induced by $P_l \xrightarrow{{\rm Id}}
P_l$, by Lemma 4 the cone is isomorphic to $P_m
\xrightarrow{\beta^{m-x_{\phi_{\Gamma}(j)}}}
P_{x_{\phi_{\Gamma}(j)}}$, by iterated application of Lemma 4 we
get: $T_{j} \simeq P_{m}$.

Application of $(\mu^-_2)^{\phi_{\Gamma}(2)} \circ \dots \circ
(\mu^-_{j-1})^{\phi_{\Gamma}(j-1)}$ gives $T_{j-1} \simeq P_{j-1},
\dots, T_{2} \simeq P_{2}$, $T_{1} \simeq P_{1}$.

Finally, $T \simeq (\mu^-_{j})^2\circ (\mu^-_{j+1})^2\circ \dots \circ
(\mu^-_{m-1})^2(A)$.\hfill\(\Box\)

Note that the indices $m$ and $j$ have switched places.

\textbf{Case 2:} $j$ is incident to the root. This case is analogous to the case
(1). Set $P_{x_{\phi_{\Gamma}(j)}}=0$ in (1). Clearly $T^{r+1}_j
\simeq P_m \xrightarrow{\beta^{m-l}} P_{l}$. By computations similar
to the case (1) we get: $\tilde{\mu}^{-(r+1)}\circ \mu^{+}_{j}\circ
\tilde{\mu}^{r}(A) \simeq (\mu^-_{j})^2\circ (\mu^-_{j+1})^2\circ
\dots \circ (\mu^-_{m-1})^2(A)$.\hfill\(\Box\)

\textbf{Case 3:} $j$ is not incident to the root and there is no edge of the same
level as $j$, following $j$ in the cyclic ordering.

$\Gamma={\Gamma^{r}}$ and ${\Gamma^{r+1}}$ are roughly of the
following form (labels $\{j+1,\dots,m-1\}$ and $\{m+1,\dots,h\}$ are on
arbitrary trees):

$$
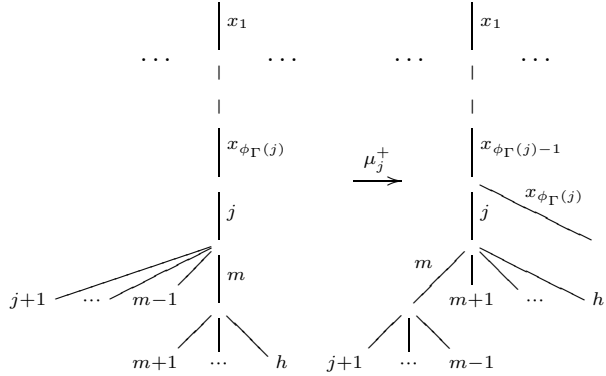
\xymatrix @! @=0.01pc {
  &  &  & \ar@{-}[d]^{x_1} &   &  &  & \ar@{-}[d]^{x_1} & & \\
&&\dots& \ar@{-}@{--}[d]  & \dots   &&\dots& \ar@{-}@{--}[d]  &\dots  & \\
  & & & \ar@{-}[d]^{x_{\phi_{\Gamma}(j)}}  &  & & & \ar@{-}[d]^{x_{\phi_{\Gamma}(j)-1}}  &&\\
  & & & \ar@{-}[d]^{j}  & & \ar[r]^{\mu^{+}_{j}} & & \ar@{-}[d]^{j} \ar@{-}[drr]^{x_{\phi_{\Gamma}(j)}}  &&\\
    & & & \ar@{-}[d]^{m}  \ar@{-}[dlll] \ar@{-}[dll] \ar@{-}[dl]   & & & & \ar@{-}[d] \ar@{-}[dl]_{m} \ar@{-}[dr] \ar@{-}[d] \ar@{-}[drr]  &&\\
     ^{j+1} & ^{\cdots} & ^{m-1} & \ar@{-}[dl]   \ar@{-}[dr] \ar@{-}[d] & & &  \ar@{-}[d] \ar@{-}[dl] \ar@{-}[dr]& ^{m+1} & ^{\cdots}&^{h}\\
   &    & ^{m+1} &  ^{\cdots} & ^{h}&^{j+1}&^{\cdots}&   ^{m-1} & & \\
}$$

As before denote by $x_1,x_2,\dots,x_{\phi_{\Gamma}(j)}=\psi(j)=j-1$
edges belonging to the shortest path from the exceptional vertex to
$j$ indexed by numbers from 1 to $\phi_{\Gamma}(j)$.

\emph{Claim:} \emph{ $\tilde{\mu}^{-(r+1)}\circ \mu^{+}_{j}\circ
\tilde{\mu}^{r}(A) \simeq F_j(A)$, where $F_j(A)$ is defined in the
Appendix (the cyclic ordering of the edges of the Brauer star
corresponding to the summands of $F_j(A)$is given by the linear
order of these summand from the bottom to the top).}

\emph{Proof.} Let us compute $T^{r+1}_j.$ The minimal left
approximation of $T^{r}_j$ is
$$\xymatrixcolsep{4pc}\xymatrix  {
  0 \ar[r] \ar[d] & P_j \ar[r]^{\beta} \ar[d]^{{\rm Id}} & P_{x_{\phi_{\Gamma}(j)}} \ar[r] \ar[d]^{soc} & 0 \ar[d] \\
  P_m \ar[r]^{\beta^{m-j}}& P_j\ar[r]^{0} & P_{x_{\phi_{\Gamma}(j)}}\ar[r]^{\beta^a} & P_{x_{\phi_{\Gamma}(j)-1}}. \\
}
$$
By lemma 4 $$T^{r+1}_j \simeq P_m
\xrightarrow{\beta^{m-x_{\phi_{\Gamma}(j)}}}
P_{x_{\phi_{\Gamma}(j)}} \xrightarrow{soc} P_{x_{\phi_{\Gamma}(j)}}
\xrightarrow{\beta^a} P_{x_{\phi_{\Gamma}(j)-1}},$$ where
$a=x_{\phi_{\Gamma}(j)}-x_{\phi_{\Gamma}(j)-1}$.

The levels of edges in $\Gamma^{r+1}$ are: $\phi_{\Gamma^{r+1}}(j) =
\phi_{\Gamma}(j)-1$, $\phi_{\Gamma^{r+1}}(j+1) =
\phi_{\Gamma}(j+1),\dots, \phi_{\Gamma^{r+1}}(m-1) =
\phi_{\Gamma}(m-1)$, $\phi_{\Gamma^{r+1}}(m) = \phi_{\Gamma}(m)-1$,
$\phi_{\Gamma^{r+1}}(m+1) = \phi_{\Gamma}(m+1)-2,\dots,
\phi_{\Gamma^{r+1}}(h) = \phi_{\Gamma}(h)-2$, hence
$$\tilde{\mu}^{-(r+1)} = (\mu^-_2)^{\phi_{\Gamma}(2)} \circ \dots \circ
(\mu^-_{j-2})^{\phi_{\Gamma}(j-2)} \circ
(\mu^-_{j})^{\phi_{\Gamma}(j)-1} \circ
(\mu^-_{m})^{\phi_{\Gamma}(m)-1} \circ $$ $$\circ
(\mu^-_{j+1})^{\phi_{\Gamma}(j+1)} \circ \dots \circ
(\mu^-_{m-1})^{\phi_{\Gamma}(m-1)} \circ
(\mu^-_{m+1})^{\phi_{\Gamma}(m+1)-2} \circ \dots \circ
(\mu^-_{h})^{\phi_{\Gamma}(h)-2} \circ $$ $$\circ
(\mu^-_{x_{\phi_{\Gamma}(j)}})^{\phi_{\Gamma}(j)} \circ
(\mu^-_{h+1})^{\phi_{\Gamma}(h+1)} \circ \dots \circ
(\mu^-_{n})^{\phi_{\Gamma}(n)}.$$ Let us apply this sequence of
mutations to $T^{r+1}$.

Applying $(\mu^-_{x_{\phi_{\Gamma}(j)}})^{\phi_{\Gamma}(j)} \circ
(\mu^-_{h+1})^{\phi_{\Gamma}(h+1)} \circ \dots \circ
(\mu^-_{n})^{\phi_{\Gamma}(n)}$ we get: $T_n \simeq P_n, \dots,
T_{h+1} \simeq P_{h+1}, T_{x_{\phi_{\Gamma}(j)}} \simeq
P_{x_{\phi_{\Gamma}(j)}}$.

Application of $(\mu^-_{h})^{\phi_{\Gamma}(h)-2}$: note that the
edge $h$ is pendant; if $h$ is not incident to $j$, then after the
application of $(\mu^-_{h})^{\phi_{\Gamma}(h)-\phi_{\Gamma}(j)-2}$,
the edge $h$ becomes incident to $j$, the complex corresponding to
that edge is $P_h \xrightarrow{\beta^{h-m}} P_m$. The minimal right
approximation of this complex is the morphism from $$P_m
\xrightarrow{\beta^{m-x_{\phi_{\Gamma}(j)}}}
P_{x_{\phi_{\Gamma}(j)}} \xrightarrow{soc} P_{x_{\phi_{\Gamma}(j)}}
\xrightarrow{\beta^a} P_{x_{\phi_{\Gamma}(j)-1}}$$ induced by $P_m
\xrightarrow{{\rm Id}} P_m$, by Lemma 4 its cone is isomorphic to
$$P_h \xrightarrow{\beta^{h-x_{\phi_{\Gamma}(j)}}}
P_{x_{\phi_{\Gamma}(j)}} \xrightarrow{soc} P_{x_{\phi_{\Gamma}(j)}}
\xrightarrow{\beta^a} P_{x_{\phi_{\Gamma}(j)-1}}.$$ The minimal
right approximation of the shift of this cone is the morphism from
$P_{x_{\phi_{\Gamma}(j)-1}} \xrightarrow{\beta^b}
P_{x_{\phi_{\Gamma}(j)-2}}$ (where $b=x_{\phi_{\Gamma}(j)-1}-
x_{\phi_{\Gamma}(j)-2}$), induced by
$P_{x_{\phi_{\Gamma}(j)-1}}~\xrightarrow{{\rm
Id}}~P_{x_{\phi_{\Gamma}(j)-1}}$, by Lemma 4 its cone is isomorphic
to
$$P_h \xrightarrow{\beta^{h-x_{\phi_{\Gamma}(j)}}}
P_{x_{\phi_{\Gamma}(j)}} \xrightarrow{soc} P_{x_{\phi_{\Gamma}(j)}}
\xrightarrow{\beta^{a+b}} P_{x_{\phi_{\Gamma}(j)-2}}.$$ By iterated
application of Lemma 4 we get: $$T_h \simeq P_h
\xrightarrow{\beta^{h-x_{\phi_{\Gamma}(j)}}}
P_{x_{\phi_{\Gamma}(j)}} \xrightarrow{soc}
P_{x_{\phi_{\Gamma}(j)}}.$$

Application of
$(\mu^-_{h-1})^{\phi_{\Gamma}(h-1)-2},\dots,(\mu^-_{m+1})^{\phi_{\Gamma}(m+1)-2}$
is completely similar and we get:
$$T_{h-1} \simeq P_{h-1} \xrightarrow{\beta^{h-1-x_{\phi_{\Gamma}(j)}}}
P_{x_{\phi_{\Gamma}(j)}} \xrightarrow{soc}
P_{x_{\phi_{\Gamma}(j)}},$$
$$\dots$$
$$T_{m+1} \simeq P_{m+1}
\xrightarrow{\beta^{m+1-x_{\phi_{\Gamma}(j)}}}
P_{x_{\phi_{\Gamma}(j)}} \xrightarrow{soc}
P_{x_{\phi_{\Gamma}(j)}}.$$

Application of $(\mu^-_{m-1})^{\phi_{\Gamma}(m-1)}$: if the edge
$m-1$ is not incident to $m$, then after the application of
$(\mu^-_{m-1})^{\phi_{\Gamma}(m-1)-\phi_{\Gamma}(m)}$, the edge
$m-1$ becomes incident to $m$, the complex corresponding to that
edge is $P_{m-1} \xrightarrow{\beta^{m-1-j}} P_j$. The minimal right
approximation of $P_{m-1} \xrightarrow{\beta^{m-1-j}} P_j$ is the
morphism $(\beta, {\rm Id})$ from $P_{m}~\xrightarrow{\beta^{m-j}}~P_j$ to $P_{m-1} \xrightarrow{\beta^{m-1-j}} P_j$, by Lemma 4 its
cone is isomorphic to $P_{m} \xrightarrow{\beta} P_{m-1}$. The
minimal right approximation of $P_{m} \xrightarrow{\beta} P_{m-1}$
is the morphism $(soc, 0,0,0)$ from $$P_m
\xrightarrow{\beta^{m-x_{\phi_{\Gamma}(j)}}}
P_{x_{\phi_{\Gamma}(j)}} \xrightarrow{soc} P_{x_{\phi_{\Gamma}(j)}}
\xrightarrow{\beta^a} P_{x_{\phi_{\Gamma}(j)-1}}$$ to $P_{m}
\xrightarrow{\beta} P_{m-1}$,
its cone is $$P_m \xrightarrow{\left( \begin{smallmatrix} -\beta^{m-x_{\phi_{\Gamma}(j)}} \\  soc\\
\end{smallmatrix} \right)} P_{x_{\phi_{\Gamma}(j)}} \oplus P_m  \xrightarrow{\left( \begin{smallmatrix} -soc & 0  \\ 0 & \beta\\
\end{smallmatrix} \right)} P_{x_{\phi_{\Gamma}(j)}} \oplus P_{m-1}
\xrightarrow{(-\beta^a,0)}P_{x_{\phi_{\Gamma}(j)-1}}.$$ The minimal
right approximation of the shift of this complex is the morphism
from $P_{x_{\phi_{\Gamma}(j)-1}} \xrightarrow{\beta^b}
P_{x_{\phi_{\Gamma}(j)-2}}$ induced by $P_{x_{\phi_{\Gamma}(j)-1}}
\xrightarrow{{\rm Id}} P_{x_{\phi_{\Gamma}(j)-1}}$, by Lemma 4 its cone is isomorphic to $$P_m \xrightarrow{\left( \begin{smallmatrix} -\beta^{m-x_{\phi_{\Gamma}(j)}} \\  soc\\
\end{smallmatrix} \right)} P_{x_{\phi_{\Gamma}(j)}} \oplus P_m  \xrightarrow{\left( \begin{smallmatrix} -soc & 0  \\ 0 & \beta\\
\end{smallmatrix} \right)} P_{x_{\phi_{\Gamma}(j)}} \oplus P_{m-1}
\xrightarrow{(-\beta^{a+b},0)}P_{x_{\phi_{\Gamma}(j)-2}}.$$ By iterated application of Lemma 4 we get: $$T_{m-1} \simeq P_m \xrightarrow{\left( \begin{smallmatrix} -\beta^{m-x_{\phi_{\Gamma}(j)}} \\  soc\\
\end{smallmatrix} \right)} P_{x_{\phi_{\Gamma}(j)}} \oplus P_m  \xrightarrow{\left( \begin{smallmatrix} -soc & 0  \\ 0 & \beta\\
\end{smallmatrix} \right)} P_{x_{\phi_{\Gamma}(j)}} \oplus P_{m-1}.$$

Similarly,
$$T_{m-2} \simeq P_m \xrightarrow{\left( \begin{smallmatrix} -\beta^{m-x_{\phi_{\Gamma}(j)}} \\  soc\\
\end{smallmatrix} \right)} P_{x_{\phi_{\Gamma}(j)}} \oplus P_m  \xrightarrow{\left( \begin{smallmatrix} -soc & 0  \\ 0 & \beta^2\\
\end{smallmatrix} \right)} P_{x_{\phi_{\Gamma}(j)}} \oplus P_{m-2},$$
$$\dots$$
$$T_{j+1} \simeq P_m \xrightarrow{\left( \begin{smallmatrix} -\beta^{m-x_{\phi_{\Gamma}(j)}} \\  soc\\
\end{smallmatrix} \right)} P_{x_{\phi_{\Gamma}(j)}} \oplus P_m  \xrightarrow{\left( \begin{smallmatrix} -soc & 0  \\ 0 & \beta^{m-j-1}\\
\end{smallmatrix} \right)} P_{x_{\phi_{\Gamma}(j)}} \oplus P_{j+1}.$$

Application of $(\mu^-_{m})^{\phi_{\Gamma}(m)-1}$: the minimal right
approximation of $P_{m}~\xrightarrow{\beta^{m-j}}~P_{j}$ is the
morphism $(soc, 0,0,0)$ from $$P_m
\xrightarrow{\beta^{m-x_{\phi_{\Gamma}(j)}}}
P_{x_{\phi_{\Gamma}(j)}} \xrightarrow{soc} P_{x_{\phi_{\Gamma}(j)}}
\xrightarrow{\beta^a} P_{x_{\phi_{\Gamma}(j)-1}}$$ to $P_{m}
\xrightarrow{\beta^{m-j}} P_{j}$, similar to the previous case we
get:
$$T_{m} \simeq P_m \xrightarrow{\left( \begin{smallmatrix} -\beta^{m-x_{\phi_{\Gamma}(j)}} \\  soc\\
\end{smallmatrix} \right)} P_{x_{\phi_{\Gamma}(j)}} \oplus P_m  \xrightarrow{\left( \begin{smallmatrix} -soc & 0  \\ 0 & \beta^{m-j}\\
\end{smallmatrix} \right)} P_{x_{\phi_{\Gamma}(j)}} \oplus P_{j}.$$

Application of $(\mu^-_{j})^{\phi_{\Gamma}(j)-1}$: the minimal right
approximation of $$P_m \xrightarrow{\beta^{m-x_{\phi_{\Gamma}(j)}}}
P_{x_{\phi_{\Gamma}(j)}} \xrightarrow{soc} P_{x_{\phi_{\Gamma}(j)}}
\xrightarrow{\beta} P_{x_{\phi_{\Gamma}(j)-1}}$$ is the morphism from
$P_{x_{\phi_{\Gamma}(j)-1}} \xrightarrow{\beta^b}
P_{x_{\phi_{\Gamma}(j)-2}}$, induced by $P_{x_{\phi_{\Gamma}(j)-1}}
\xrightarrow{{\rm Id}} P_{x_{\phi_{\Gamma}(j)-1}}$, by iterated
application of Lemma 4 we get: $$T_{j} \simeq P_m
\xrightarrow{\beta^{m-x_{\phi_{\Gamma}(j)}}}
P_{x_{\phi_{\Gamma}(j)}} \xrightarrow{soc}
P_{x_{\phi_{\Gamma}(j)}}.$$

Application of $(\mu^-_2)^{\phi_{\Gamma}(2)} \circ \dots \circ
(\mu^-_{j-2})^{\phi_{\Gamma}(j-2)}$ clearly gives $$T_{j-2} \simeq
P_{j-2},\dots, T_{2} \simeq P_2, T_{1} \simeq P_1.$$\hfill\(\Box\)

\emph{Claim:} \emph{ $F_j(A)=(\mu^+_{h})^2 \circ \dots \circ
(\mu^+_{m+1})^2 \circ (\mu^-_{j})^2 \circ (\mu^-_{j+1})^2 \circ \dots
\circ (\mu^-_{m-1})^2 \circ (\mu^+_{m})^2 \circ (\mu^+_{m-1})^2
\circ \dots \circ (\mu^+_{j})^2(A)$.}

\emph{Proof.} Let us apply to $F_j(A)$ the inverse sequence of
mutations. If we obtain $A$, the claim is proved. Let us compute a
double mutation of $T_h$ along $T_{x_{\phi_{\Gamma}(j)}}$. The
minimal right approximation of $T_h$ is the morphism from
$T_{x_{\phi_{\Gamma}(j)}}$ to $T_h$ induced by
$P_{x_{\phi_{\Gamma}(j)}} \xrightarrow{{\rm Id}}
P_{x_{\phi_{\Gamma}(j)}}$, by Lemma 4 its cone is isomorphic to $P_h
\xrightarrow{\beta^{h-x_{\phi_{\Gamma}(j)}}}
P_{x_{\phi_{\Gamma}(j)}}$; after applying the shift
$P_{x_{\phi_{\Gamma}(j)}}$ is concentrated in the same degree as in
$T_{x_{\phi_{\Gamma}(j)}}$, applying Lemma 4 one more time we get:
$$[(\mu^-_{h})^2(F_j(A))]_h \simeq P_h.$$

Similarly,
$$[(\mu^-_{h-1})^2 \circ (\mu^-_{h})^2(F_j(A))]_{h-1} \simeq
P_{h-1}$$
$$\dots$$
$$[(\mu^-_{m+1})^2 \circ\dots\circ (\mu^-_{h-1})^2
\circ (\mu^-_{h})^2(F_j(A))]_{m+1} \simeq P_{m+1}.$$

Let us denote the complex $(\mu^-_{m+1})^2 \circ\dots\circ
(\mu^-_{h-1})^2 \circ (\mu^-_{h})^2(F_j(A))$ by $T^r$.

Recall that the homotopy category can be defined as the stable
category of the Frobenius category of complexes with respect to
degree-wise split short exact sequences. So the distinguished
triangles in the homotopy category come from degree-wise split short
exact sequences \cite{Ha}.

Let us compute a double mutation of $T_m$ along $T_{j}$. Consider
the following degree-wise split short exact sequence in the category
of complexes: {\small $$
\xymatrixcolsep{2pc}\xymatrixrowsep{4pc}\xymatrix {
X= \ar[d]^{g} &  0 \ar[rr] \ar[d] && P_m \ar[rr]^{\beta^{m-j}} \ar[d]^{\left( \begin{smallmatrix} 0 \\  {\rm Id}\\
\end{smallmatrix} \right)} && P_j \ar[d]^{\left( \begin{smallmatrix} 0 \\  {\rm Id}\\
\end{smallmatrix} \right)} \\
T_m= \ar[d]^{f}  &   P_m \ar[d]^{{\rm Id}} \ar[rr]^(0.4){\left( \begin{smallmatrix} -\beta^{m-x_{\phi_{\Gamma}(j)}} \\  soc\\
\end{smallmatrix} \right)} && P_{x_{\phi_{\Gamma}(j)}} \oplus P_m \ar[d]^{({\rm Id},0)} \ar[rr]^{\left( \begin{smallmatrix} -soc & 0  \\ 0 & \beta^{m-j}\\
\end{smallmatrix} \right)}&& P_{x_{\phi_{\Gamma}(j)}} \oplus P_j \ar[d]^{({\rm Id},0)} \\
T_j \simeq &   P_m \ar[rr]^{-\beta^{m-x_{\phi_{\Gamma}(j)}}} && P_{x_{\phi_{\Gamma}(j)}} \ar[rr]^{-soc} && P_{x_{\phi_{\Gamma}(j)}},\\
}$$} where $f$ is clearly a minimal left approximation of $T_m$ with
respect to other summands of $T^r$. So $Cone(f) \simeq (P_m
\xrightarrow{\beta^{m-j}} P_j)$ concentrated in the corresponding
degree. The minimal left approximation of $Cone(f)$ with respect to
other summands of $T^{r+1}$ is the morphism to $$T_j = P_m
\xrightarrow{\beta^{m-x_{\phi_{\Gamma}(j)}}}
P_{x_{\phi_{\Gamma}(j)}} \xrightarrow{soc} P_{x_{\phi_{\Gamma}(j)}}$$
induced by $({\rm Id},\beta^{j-x_{\phi_{\Gamma}(j)}})$, by Lemma 4
its cone is isomorphic to $$[(\mu^+_{m})^2(T^r)]_m \simeq P_j
\xrightarrow{\beta^{j-x_{\phi_{\Gamma}(j)}}}
P_{x_{\phi_{\Gamma}(j)}} \xrightarrow{soc}
P_{x_{\phi_{\Gamma}(j)}}.$$

Similarly,
$$[(\mu^+_{j+1})^2
\circ (\mu^+_{m})^2(T^r)]_{j+1} \simeq P_{j+1}
\xrightarrow{\beta^{j+1-x_{\phi_{\Gamma}(j)}}}
P_{x_{\phi_{\Gamma}(j)}} \xrightarrow{soc}
P_{x_{\phi_{\Gamma}(j)}}$$
$$\dots$$
$$[(\mu^+_{m-1})^2 \circ\dots\circ (\mu^+_{j+1})^2
\circ (\mu^+_{m})^2(T^r))]_{m-1} \simeq P_{m-1}
\xrightarrow{\beta^{m-1-x_{\phi_{\Gamma}(j)}}}
P_{x_{\phi_{\Gamma}(j)}} \xrightarrow{soc}
P_{x_{\phi_{\Gamma}(j)}}.$$

So we have:
$$(\mu^+_{m-1})^2 \circ\dots\circ (\mu^+_{j+1})^2
\circ (\mu^+_{m})^2 \circ (\mu^-_{m+1})^2 \circ\dots\circ
(\mu^-_{h-1})^2 \circ (\mu^-_{h})^2(F_j(A)) \simeq G_j(A),$$ where
$G_j(A)$ is defined in the Appendix (the cyclic ordering of the
edges corresponding to the summands of $G_j(A)$ in the Brauer star
is defined by the linear order of these summands from the bottom to
the top).

Applying $(\mu^-_{j})^2$ to $G_j$ we get $[(\mu^-_{j})^2(G_j)]_j
\simeq P_m$.

Similarly,
$$[(\mu^-_{m-1})^2 \circ (\mu^-_{j})^2(G_j)]_{m-1} \simeq P_{m-1}$$
$$\dots$$
$$[(\mu^-_{j+1})^2 \circ\dots(\mu^-_{m-1})^2 \circ (\mu^-_{j})^2(G_j)]_{j+1} \simeq P_{j+1}$$
$$[(\mu^-_{m})^2 \circ (\mu^-_{j+1})^2 \circ\dots(\mu^-_{m-1})^2 \circ (\mu^-_{j})^2(G_j)]_{m} \simeq P_j.$$
So $(\mu^-_{m})^2 \circ (\mu^-_{j+1})^2 \circ\dots(\mu^-_{m-1})^2
\circ (\mu^-_{j})^2(G_j) \simeq A$.

Note that in the resulting copy of $A$ the cyclic ordering of the
projective modules is standard, but the labelling is not: $j$ and
$m$ have switched their places. Hence the assertion
follows.\hfill\(\Box\)

This finishes the proof of Theorem 3.

\subsection{Multiplicity free case}

Let $A$ be a Brauer star algebra of type $(n,1)$. Consider a
subgroup $\tilde{\mathcal{R}}$ of the derived Picard group of $A$
generated by shift, ${\rm Pic}(A)$, equivalences $H_i$ and
equivalences $Q_i$, where $Q_i$ are some standard equivalences
acting on the indecomposable projective modules as follows:
$$Q_i(P_j)=\left\{
            \begin{array}{ll}
              0\mbox{ } \rightarrow \mbox{ } P_i\mbox{ } \mbox{ }\rightarrow \mbox{ } 0, & j=i  \\
              0\mbox{ }\rightarrow \mbox{ } P_i\mbox{ } \mbox{ }\xrightarrow{\beta^{i-j}} \mbox{ } P_{j}, & j \neq i, \\
            \end{array}
          \right.$$
where $P_i$ are concentrated in degree $0$.

\begin{rem}
The group generated by $Q_i$ was considered by Muchtadi-Alamsyah in
\cite{IM}. It has an action of the braid group whose associated
Coxeter group is given by the complete graph on $n$ vertices, but
this action is not faithful.
\end{rem}
\begin{rem}
It is clear that $Q_i(A) \simeq \mu^-_{i+1} \circ \mu^-_{i+2} \circ
\dots \circ \mu^-_{i-2} \circ \mu^-_{i-1}(A)$, this series of
mutations changes the central vertex of the Brauer star to the outer
vertex of the edge labelled $i$. It is also clear that one can
obtain $Q_i$ from $Q_{i+1}$ by conjugation with the automorphism
corresponding to the rotation of the Brauer star. Note also that
$Q_i^2$ does not change the central vertex of the Brauer star, hence
$Q_i^2 \in {\mathcal{R}}$.
\end{rem}

\begin{thm}
If $t=1$, then ${\rm TrPic}(A) = \tilde{\mathcal{R}}$.
\end{thm}

\textbf{Proof.} Since $Q_i(A) \simeq \mu^-_{i+1} \circ \mu^-_{i+2}
\circ \dots \circ \mu^-_{i-2} \circ \mu^-_{i-1}(A)$, the summands of
$Q_i(A)$ can be computed as cones of minimal approximations of
projective modules. Similar to the computation of $F(H_i(A))$, in
order to compute $F(Q_i(A))$ for some autoequivalence $F$ one can
apply the sequence of mutations $\mu^-_{i+1} \circ \mu^-_{i+2} \circ
\dots \circ \mu^-_{i-2} \circ \mu^-_{i-1}$ to $F(A)$.

As before, it is sufficient to prove that for any tilting complex
$T$ such that ${\rm End}_{D^b(A)}(T)^{op} \simeq A$ there is an
element from $\tilde{\mathcal{R}}$ which sends the indecomposable
projective $A$-modules to the summands of $T$. Assume that $T$ is
concentrated in non-positive degrees. By results of
Aihara~\cite{Ai2} $T=\mu^{+}_{i_s} \circ\dots\circ \mu^{+}_{i_{r+1}}
\circ \mu^{+}_{i_r} \circ \mu^{+}_{i_{r-1}} \circ\dots \circ
\mu^{+}_{i_2} \circ \mu^{+}_{i_1}(A)$ for some $(i_1, i_2, \dots,
i_s)$. If the series of mutations $\mu^{+}_{i_s}
 \circ\dots\circ \mu^{+}_{i_{r+1}} \circ  \mu^{+}_{i_r} \circ
 \mu^{+}_{i_{r-1}} \circ\dots \circ \mu^{+}_{i_2} \circ
 \mu^{+}_{i_1}(A)$ does not change the central vertex of the Brauer
 star, then by calculations of the previous section the standard
 equivalence $F$ such that $F(A) \simeq T$
 belongs to $\mathcal{R}$ and hence to $\tilde{\mathcal{R}}$. If the series
 of mutations changes the central vertex of the Brauer
 star to the outer vertex of the edge labelled $i$, then applying
 $\mu^-_{i+1} \circ \mu^-_{i+2} \circ \dots \circ \mu^-_{i-2} \circ
\mu^-_{i-1}$ we obtain a tilting complex $T'$ such that there is some $F' \in
 \mathcal{R},$ $F'(A)\simeq T'$. So we have $T' \simeq  F \circ Q_i (A)$, hence $F$ belongs to $\tilde{\mathcal{R}}$.\hfill\(\Box\)

\begin{rem}
As in the case $t>1$, the question of the relations in
$\tilde{\mathcal{R}}$ remains open.
\end{rem}

\begin{rem}
It seems that the action of the derived Picard group on the edges of
the Brauer star can not be defined properly, since in many examples
of computations we have found two sequences of mutations which give
the same tilting complex but act differently on the edges of the
Brauer star.
\end{rem}

\section{Appendix}
$$F_j= \left( \begin{smallmatrix}
                                           T_1 \\
                                           T_2 \\
                                           \vdots \\
                                           T_{j-2} \\
                                           T_j \\
                                           T_m \\
                                           T_{j+1} \\
                                           \vdots \\
                                           T_{m-1} \\
                                           T_{m+1} \\
                                           \vdots \\
                                           T_h \\
                                           T_{x_{\phi_{\Gamma}(j)}} \\
                                           T_{h+1} \\
                                           \vdots \\
                                           T_n
                                         \end{smallmatrix} \right) =
\left( \begin{smallmatrix}
      &   &   &   & P_1 \\
      &   &   &   & P_2 \\
      &   &   &   & \vdots \\
      &   &   &   & P_{j-2} \\
    P_m & \xrightarrow{\beta^{m-x_{\phi_{\Gamma}(j)}}} & P_{x_{\phi_{\Gamma}(j)}} & \xrightarrow{soc} & P_{x_{\phi_{\Gamma}(j)}} \\
    P_m & \xrightarrow{\left( \begin{smallmatrix} -\beta^{m-x_{\phi_{\Gamma}(j)}} \\  soc\\
\end{smallmatrix} \right)} & P_{x_{\phi_{\Gamma}(j)}} \oplus P_m & \xrightarrow{\left( \begin{smallmatrix} -soc & 0  \\ 0 & \beta^{m-j}\\
\end{smallmatrix} \right)} & P_{x_{\phi_{\Gamma}(j)}} \oplus P_j \\
    P_m & \xrightarrow{\left( \begin{smallmatrix} -\beta^{m-x_{\phi_{\Gamma}(j)}} \\  soc\\
\end{smallmatrix} \right)} & P_{x_{\phi_{\Gamma}(j)}} \oplus P_m & \xrightarrow{\left( \begin{smallmatrix} -soc & 0  \\ 0 & \beta^{m-j-1}\\
\end{smallmatrix} \right)} & P_{x_{\phi_{\Gamma}(j)}} \oplus P_{j+1} \\
      &   &   &   & \vdots \\
    P_m & \xrightarrow{\left( \begin{smallmatrix} -\beta^{m-x_{\phi_{\Gamma}(j)}} \\  soc\\
\end{smallmatrix} \right)} & P_{x_{\phi_{\Gamma}(j)}} \oplus P_m & \xrightarrow{\left( \begin{smallmatrix} -soc & 0  \\ 0 & \beta\\
\end{smallmatrix} \right)} & P_{x_{\phi_{\Gamma}(j)}} \oplus P_{m-1} \\
    P_{m+1} & \xrightarrow{\beta^{m+1-x_{\phi_{\Gamma}(j)}}} & P_{x_{\phi_{\Gamma}(j)}} & \xrightarrow{soc} & P_{x_{\phi_{\Gamma}(j)}} \\
      &   &   &   & \vdots \\
    P_{h} & \xrightarrow{\beta^{h-x_{\phi_{\Gamma}(j)}}} & P_{x_{\phi_{\Gamma}(j)}} & \xrightarrow{soc} & P_{x_{\phi_{\Gamma}(j)}} \\
      &   &   &   & P_{x_{\phi_{\Gamma}(j)}} \\
      &   &   &   & P_{h+1} \\
      &   &   &   & \vdots \\
      &   &   &   & P_n \\
  \end{smallmatrix}
\right)
                                   $$

$$G_j= \left( \begin{smallmatrix}
                                           T_1 \\
                                           T_2 \\
                                           \vdots \\
                                           T_{j-2} \\
                                           T_m \\
                                           T_{j+1} \\
                                           \vdots \\
                                           T_{m-1} \\
                                            T_j \\
                                             T_{x_{\phi_{\Gamma}(j)}} \\
                                           T_{m+1} \\
                                           \vdots \\
                                           T_h \\
                                           T_{h+1} \\
                                           \vdots \\
                                           T_n
                                         \end{smallmatrix} \right) =
\left( \begin{smallmatrix}
      &   &   &   & P_1 \\
      &   &   &   & P_2 \\
      &   &   &   & \vdots \\
      &   &   &   & P_{j-2} \\
    P_j & \xrightarrow{\beta^{j-x_{\phi_{\Gamma}(j)}}} & P_{x_{\phi_{\Gamma}(j)}} & \xrightarrow{soc} & P_{x_{\phi_{\Gamma}(j)}} \\
    P_{j+1} & \xrightarrow{\beta^{j+1-x_{\phi_{\Gamma}(j)}}} & P_{x_{\phi_{\Gamma}(j)}}  & \xrightarrow{soc} & P_{x_{\phi_{\Gamma}(j)}}\\
      &   &   &   & \vdots \\
    P_{m-1} & \xrightarrow{\beta^{m-1-x_{\phi_{\Gamma}(j)}}} & P_{x_{\phi_{\Gamma}(j)}}  & \xrightarrow{soc} & P_{x_{\phi_{\Gamma}(j)}}\\
    P_{m} & \xrightarrow{\beta^{m-x_{\phi_{\Gamma}(j)}}} & P_{x_{\phi_{\Gamma}(j)}} & \xrightarrow{soc} & P_{x_{\phi_{\Gamma}(j)}} \\
 &   &   &   & P_{x_{\phi_{\Gamma}(j)}} \\
 &   &   &   & P_{m+1} \\
      &   &   &   & \vdots \\
   &   &   &   & P_{h} \\
      &   &   &   & P_{h+1} \\
      &   &   &   & \vdots \\
      &   &   &   & P_n \\
  \end{smallmatrix}
\right)
                                   $$

\end{document}